\documentclass[11pt, reqno]{amsart}
\usepackage[pdftex,hidelinks]{hyperref}
\usepackage[numbers]{natbib}
\usepackage{graphicx}
\usepackage[shortlabels]{enumitem}
\usepackage{fullpage}
\usepackage{enumitem}
\usepackage{xcolor}

\usepackage{todonotes}

\usepackage[ansinew]{inputenc} 

\usepackage{amssymb}
\usepackage{bbm}
\usepackage{bm} 

\usepackage{mathtools}
\usepackage{amsthm}
\usepackage{cases}
\usepackage{MnSymbol}

\newcommand{\Gau}{\mathfrak{G}} 
\newcommand{\Sti}{G} 
\newcommand{\U}{U} 
\newcommand{\intprocleft}{\mathfrak{F}_{1}}
\newcommand{\intprocleftup}{\mathfrak{F}}

\renewcommand{\Re}{\operatorname{Re}}

\newcommand{\E}{\mathbb E }
\newcommand{\R}{\mathbb{R}}
\newcommand{\N}{\mathbb{N}}
\newcommand{\C}{\mathbb{C}}

\newcommand{\1}{\mathbbm{1}}

\newcommand{\Cov}{\mathop{\mathrm{Cov}}\nolimits}

\newcommand{\eee}{{\rm e}}

\newcommand{\toweak}{\overset{w}{\underset{n\to\infty}\longrightarrow}}
\newcommand{\toweakk}{\overset{w}{\underset{n\to\infty}\longrightarrow}}

\newcommand{\ton}{\overset{}{\underset{n\to\infty}\longrightarrow}}

\newcommand{\dd}{{\rm d}}

\newcommand{\eps}{\varepsilon}

\theoremstyle{plain}
\newtheorem{theorem}{Theorem}[section]
\newtheorem*{theorem*}{Theorem}
\newtheorem{appl}{Application}
\newtheorem{proposition}[theorem]{Proposition}

\newtheorem{lemma}[theorem]{Lemma}
\newtheorem{corollary}[theorem]{Corollary}

\theoremstyle{definition}
\newtheorem{definition}[theorem]{Definition}

\newtheorem{example}[theorem]{Example}

\newtheorem{assu}{Assumption}
\newenvironment{ass}[1]
{%
	\renewcommand{\theassu}{#1}%
	\begin{assu}%
	}
	{%
	\end{assu}%
}

\theoremstyle{remark}
\newtheorem{remark}[theorem]{Remark}

\begin{document}

\title[Zeros and exponential profiles III: Fluctuations]{Zeros and exponential profiles\\of (random) polynomials III: Fluctuations 
}

\begin{abstract}
The first two parts of this series [Jalowy, Kabluchko and Marynych, Constructive Approximation, 2026] showed how the limiting zero distribution of a polynomial determines the exponential profile of its coefficients and vice versa. This third part studies fluctuations about that profile across coefficient indices $k=\lfloor \alpha n\rfloor$ proportional to the degree $n$. Given nonnegative random variables $(\lambda_{j:n})_{1\le j\le n}$ representing reflected zeros, we show that functional limit theorems for their empirical measures yield a functional limit theorem in the Skorokhod space for the process $\alpha\mapsto\log e_{\lfloor n\alpha\rfloor:n}$, where $e_{k:n}$ denotes the $k$-th elementary symmetric polynomial in $\lambda_{j:n}$. 

For independent and identically distributed $\lambda_{j:n}$, our result extends classical scalar central limit theorems to process-level fluctuations and further generalizes to dependent and non-identically distributed arrays. Another application concerns the eigenvalues $\lambda_{j:n}$ of sample covariance matrices. We show that the normalized coefficients of their characteristic polynomials converge to a lognormal process whose logarithm is the restriction to the real axis of a Gaussian analytic function with covariance kernel $(\alpha,\beta)\mapsto-\log(1-c\alpha\beta)$, where $c$ is the limiting ratio of the number of rows to the number of columns. This parallels previously established results on fluctuations of logarithms of characteristic polynomials.
\end{abstract}

	\author{Jonas Jalowy}
	\address{Jonas Jalowy: Paderborn University, Institute of Mathematics, Warburger Str. 100, 33098 Paderborn, Germany}
	\email{jjalowy@math.upb.de}
		
	\author{Alexander Marynych}
	\address{Alexander Marynych: Department of Mathematics, Kyiv School of Economics, Kyiv, Ukraine}
\email{omarynych@kse.org.ua}

\keywords{Legendre transform; Stieltjes transform; symmetric polynomials; Wishart characteristic polynomial}
\subjclass[2020]{Primary: 60F17, 60B10;  Secondary: 12D05, 05E05}
\maketitle

\section{Introduction}

Let $(\lambda_{j:n})_{n\in \N, 1\leq j\leq n}$ be a collection of random variables taking values in $[0,\infty)$ and defined on some probability space $(\Omega,\mathcal{F},\mathbb{P})$. We are interested in the asymptotic behavior of the elementary symmetric polynomials
\begin{equation}\label{eq:e_n_def}
	e_{k:n}:=
	\sum_{1\leq i_1<i_2<\cdots<i_k\leq n}\lambda_{i_1:n}\cdot\lambda_{i_2:n}\cdot\ldots\cdot\lambda_{i_k:n},\quad k\in\{0,1,\ldots,n\},
\end{equation}
where we stipulate that $e_{0:n}=1$. By Vieta's theorem, they are the decreasing-order coefficients of the monic polynomial with roots $-\lambda_{j:n}$ given by
\begin{align}\label{eq:p_n}
	p_n(z):=\prod_{j=1}^{n}(z+\lambda_{j:n})= \sum_{k=0}^{n}e_{n-k:n}z^k ,\qquad z\in\mathbb{C},
\end{align}
that is, $[z^{n-k}]\left(p_n(z)\right)=e_{k:n}$.

Let us present our main results and its applications in a simplified form and setting, before introducing the technicalities rigorously in the upcoming section.
Suppose that the empirical law $\mu_n$ of the $\lambda_{j:n}$ converges almost surely to a deterministic probability law $\mu_\infty$ on $(0,\infty)$. Then, the \emph{normalized logarithmic potential} converges, that is, as $n\to\infty$,
$$
\Psi_n(z):=n^{-1}\log\left(\frac{p_n(z)}{p_n(1)}\right)=\int_{[0,+\infty)}\log\left(\frac{z+x}{1+x}\right)\mu_n({\rm d}x)\ \to \int_{[0,+\infty)}\log\left(\frac{z+x}{1+x}\right)\mu_{\infty}({\rm d}x)=:\Psi_\infty(z)
$$ 
locally uniformly on the right half-plane. In Part I, see~\cite[Theorem 2.4]{JalKabMar:I}, of this series of papers, we have shown that weak convergence of the empirical measures of zeros of $p_n$ to a probability measure $\mu_{\infty}$ on $(0,+\infty)$ is equivalent to the logarithmic asymptotics of its coefficients
$$\frac1n\log\left(\frac{e_{\lfloor \alpha n\rfloor:n}}{p_n(1)}\right)\to g(1-\alpha)=\inf_{x>0}\bigl\{\Psi_\infty(x)-(1-\alpha)\log x\bigr\},\qquad \alpha\in(0,1).$$
We call $g$ the \emph{exponential profile}. Our main result of the present paper shows how the logarithmic potential transfers its fluctuations to the coefficients about their limit $g$.

\begin{theorem*}[Simplified Theorem \ref{thm:main}] 
If there is a rescaling $c_n=o(n/\log n)$ with $c_n\to\infty$ such that we have weak convergence $$c_n(\Psi_n-\Psi_\infty)\toweak\intprocleft$$ to some analytic process $\intprocleft$  in the locally uniform topology on the right half-plane, then we have locally uniform weak convergence of the process
$$
	c_n\left(\frac1n\log\left(\frac{e_{\lfloor n\alpha\rfloor:n}}{p_n(1)}\right)
	-g(1-\alpha)\right)_{\alpha\in(0,1)}\toweak\intprocleft(\eee^{-g'(1-\alpha)})_{\alpha\in(0,1)}.$$
\end{theorem*}
Since $\Psi_n$ contains all information on the zeros, this reveals once more the connection that has been established in \cite{JalKabMar:I} between the zeros and coefficients, even on the level of fluctuations.
In the upcoming Section~\ref{sec:results}, we will rigorously introduce our assumptions \ref{ass:A1} and \ref{ass:A2}. Then, we present the full theorem, which allows for mass of $\mu_\infty$ at $0$ and $\infty$, and for an $n$-dependent centering, including another deterministic term for scales $c_n\sim n$,
 which is relevant for random matrix theory.

Our probabilistic proof relies on representation \eqref{eq:coeff} below, relating the coefficients to the number of successes in $n$ inhomogeneous (conditionally) independent Bernoulli trials with random success probabilities $(1+\lambda_{j:n})^{-1}$. 
We refine the conditional large-deviation analysis of Part I \cite{JalKabMar:I} to show how the functional limit theorem for $\Psi_n$ translates into joint fluctuations of the coefficients. 
The proof will be given in Section \ref{sec:proof}. To illustrate the scope of the theorem, Section~\ref{sec:appl} presents three applications to random polynomial models. For the moment, we deliberately omit certain technical details, such as the function spaces in which convergence takes place.
 
\begin{appl}[Corollary \ref{cor:iid}]
	 For independent copies $\lambda_{j:n}$ of $\lambda$ distributed according to $\mu_\infty$ with a finite logarithmic moment of order $2+\delta$, we have convergence of the process
	 \begin{align*}
	 	\sqrt n	\left( \tfrac 1 n \log (e_{\lfloor \alpha n \rfloor:n}) - g(1-\alpha)-\E[\log(1+\lambda)]\right) &\toweak~\Gau(\alpha)
	 \end{align*}
	where $\Gau$ is a centered Gaussian process with  
	 \begin{align*}\Cov(\Gau(\alpha),\Gau(\beta))=\Cov\big(\log(\eee^{-g'(1-\alpha)}+\lambda),\log(\eee^{-g'(1-\beta)}+\lambda)\big).\end{align*} 
\end{appl}

Compared to earlier results such as \cite{Bochi+Iommi+Ponce:2022,major_symm_iid,Szekely:1982,vanEs:1986}, the above weakens integrability requirements and strengthens the scalar CLT's to functional limit theorems with explicit covariance structure, see Remark~\ref{rem:literature_iid} for bibliographic comments. Further, we generalize this to independent, non-identically distributed triangular arrays in Theorem~\ref{thm:independent}, such as deterministic roots perturbed by an i.i.d.~sequence in Corollary \ref{cor:iid-perturbations}. 

\begin{appl}[Corollary \ref{cor:RMT}]
Let ${\bf X}$ be an $n\times m$ random matrix with complex entries being i.i.d.~copies of $X$ with their first four moments matching Gaussians and ${\bf X}{\bf X}^*/m$  be the covariance matrix with $n/m\to c\in(0,\infty)$. For the decreasing coefficients $e_{k:n}$ of its characteristic polynomial, we have the convergence of processes
		$$
\frac{e_{\lfloor \alpha n\rfloor:n}}{\mathbb E e_{\lfloor\alpha n\rfloor:n}}\toweak \sqrt{1-c\alpha^2}e^{\Gau(\alpha)},\qquad \text{ on	}(0,\min\{1,1/c\})	$$
 where $\Gau(\alpha)\overset{\mathcal D}=\sum_{k\ge1} {Z_k} \sqrt{c^k/k}\alpha^k$ for i.i.d.~$Z_k\sim\mathcal N(0,1)$ is a centered Gaussian process with $${\rm Cov}(\Gau(\alpha),\Gau(\beta))=-\log(1-c\alpha\beta).$$ 
\end{appl}
We will discuss related literature in Remark \ref{rem:literature-RMT}, where we relate this result to characteristic polynomials of Gaussian and circular-$\beta$-ensembles previously discovered in \cite{LambertPaquette,Naj}, see also \cite{dumitriupaquette,HLN08,johnstone2020testing}.  Another particularly simple application will be the following.
\begin{appl}[Equation \eqref{eq:indepcoeff}]
For every $m\in\mathbb N$ and any prescribed distinct points $\alpha_1,\ldots,\alpha_m\in(0,1)$, there exists a sequence of real-rooted random polynomials whose logarithmic coefficient fluctuations converge in distribution to a Gaussian analytic process $\Gau$ such that $\Gau(\alpha_1),\dots,\Gau(\alpha_m)$ are independent standard normal random variables. In particular,
$$
c_n\left(
\frac1n\log\frac{e_{n-\lfloor\alpha_jn\rfloor:n}}{p_n(1)}
-g(\alpha_j)
\right)_{j=1}^m
~\toweak~(\Gau(\alpha_j))_{j=1}^m\sim\mathcal N_m(0,I_m).
$$
	\end{appl}
Note that standard models of independent coefficients such as \cite{Kac,LO,KZ14}  typically have complex zeros. Although the limiting Gaussian process $\Gau$ can have independent values at any prescribed finite collection of points, its analyticity excludes a nondegenerate white noise limit on an interval.

\begin{figure}[ht]
	\centering
	\includegraphics[width=0.325\columnwidth]{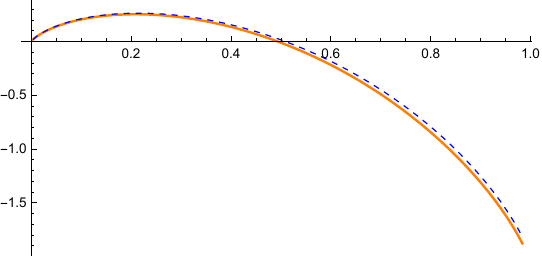}
	\includegraphics[width=0.325\columnwidth]{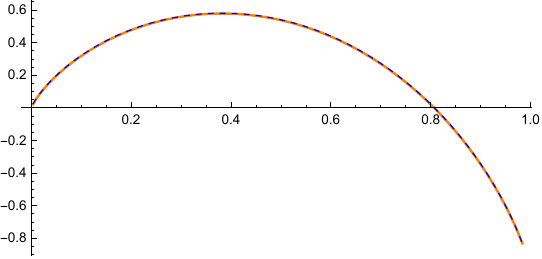}
	\includegraphics[width=0.325\columnwidth]{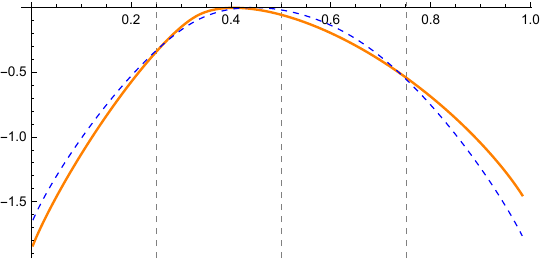}\\	\includegraphics[width=0.325\columnwidth]{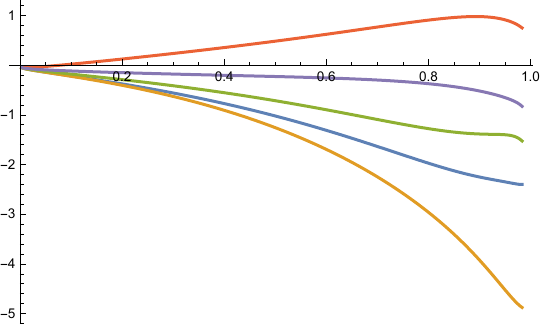}
	\includegraphics[width=0.325\columnwidth]{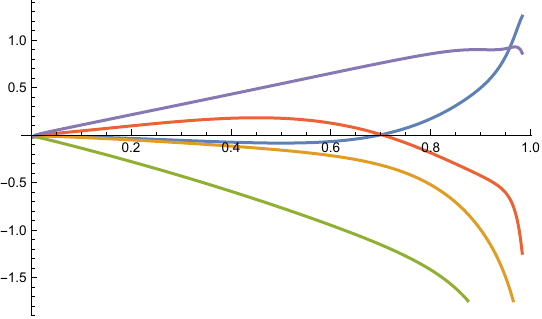}
	\includegraphics[width=0.325\columnwidth]{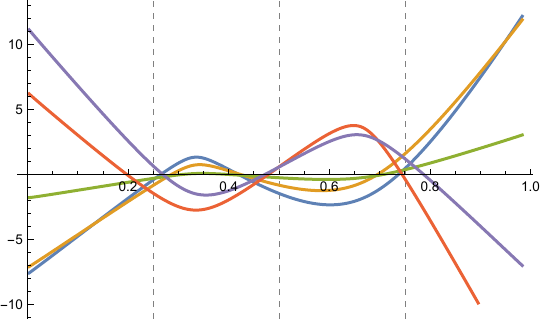}
	\caption{Simulations of the Applications 1-3 from left to right for $n=1500$. The top row shows the deterministic coefficient profile $g$ (dashed blue) and a finite $n$ approximation (orange). The bottom row illustrates the analytic limit $\Gau$ through five finite-$n$ realizations of the rescaled fluctuations about $g$. The models use $\lambda=U^2$ with $ U\sim\mathrm{Unif}(0,1)$, square complex Wishart matrices with standard Gaussian entries, and three locations $\{\alpha_1,\alpha_2,\alpha_3\}=\{1/4,2/4,3/4\}$. The first two limits are pinned at $\alpha=0$, the third is pinned at $\alpha=4/9$ due to its normalization. Observe that the fluctuation order  $1/n$ in the Wishart Application 2 shows faster convergence to $g$. As $\alpha\uparrow1$, the Wishart variance diverges and the i.i.d. variance remains bounded. In Application 3, the limiting values at $\alpha_j$ are independent standard Gaussians.}
	\label{fig:sim}
\end{figure}

\section{Main results}\label{sec:results}
\subsection{Preliminaries}

In order to cover arbitrary distributions $\mu_{\infty}$ on $[0,+\infty]$, we will need to normalize $p_n$ in \eqref{eq:p_n} as follows. Naturally, this leads to a probabilistic interpretation, which we shall need for our proofs, but not for the presentation of our results. For $n\in\mathbb{N}$, let $\mathcal{F}_n$ be the sigma-algebra generated by $\lambda_{1:n},\lambda_{2:n},\ldots,\lambda_{n:n}$. Let $B_{1:n},B_{2:n},\ldots,B_{n:n}$ be conditionally independent, given $\mathcal{F}_n$, Bernoulli random variables such that
\begin{align}\label{eq:eta_j}
	\mathbb{P}[B_{j:n}=1 | \mathcal{F}_n]=1-\mathbb{P}[B_{j:n}=0 | \mathcal{F}_n]=\eta_{j:n} :=\frac 1 {1+\lambda_{j:n} }\in (0,1],\quad 1\leq j\leq n.
\end{align}
Put $S_n:=\sum_{j=1}^{n}B_{j:n}$ for $n\in\mathbb{N}$ and define the normalized (random) polynomial $P_n$, given by the probability generating function of $S_n$
\begin{align}\label{eq:P_n}
	P_n(z):=\mathbb{E}[z^{S_n}|\mathcal{F}_n]=\prod_{j=1}^{n}\frac{z+\lambda_{j:n}}{1+\lambda_{j:n}},\qquad z\in\mathbb{C}.
\end{align} 
It has coefficients
\begin{equation}\label{eq:coeff}
	a_{k:n}:=[z^{k}]\left(P_n(z)\right)=\mathbb{P}[S_n=k | \mathcal{F}_n]=\frac{e_{n-k:n}}{p_n(1)},\quad k\in \{0,1,\ldots,n\}.
\end{equation}

Our main focus is put on the c\`{a}dl\`{a}g random process $g_n:=(g_n(\alpha))_{\alpha\in [0,1]}$ given via
\begin{equation}\label{eq:g_n_def}
	g_n(\alpha):=\frac{1}{n}\log a_{\lfloor  \alpha n\rfloor:n},\quad \alpha\in [0,1],
\end{equation}
with the convention $\log 0 = 0$. 

\subsection{Our assumptions}
We shall formulate our assumptions on $(\lambda_{j:n})_{n\in \N, 1\leq j\leq n}$ in terms of various objects such as empirical measures, Cauchy-Stieltjes transforms and logarithmic potential, and we follow the terminology and notation of \cite{JalKabMar:I,JalKabMarII}, summarized in Remark \ref{rem:notation}.

\begin{definition}
	For each $n\in\mathbb{N}$, let 
$$
\mu_n:=\frac{1}{n}\sum_{k=1}^{n}\delta_{\lambda_{k:n}}
$$
be the \emph{empirical distribution} of $\lambda_{1:n},\lambda_{2:n},\ldots,\lambda_{n:n}$. 
\end{definition}
Our first standing assumption is a.s.~weak convergence of $\mu_n$.

\begin{ass}{$({\bf A}_1)$}\label{ass:A1}
	For some deterministic  distribution $\mu_{\infty}$ on $[0,+\infty]$ satisfying $\mu_{\infty}((0,+\infty))>0$,
\begin{equation}\label{eq:A1}
\mu_n~\toweak~\mu_{\infty}\qquad\text{a.s.},
\end{equation}
on the space of probability measures on $[0,+\infty]$ endowed with the topology of weak convergence.
\end{ass}
Note that Assumption \ref{ass:A1} allows an escape of {\it some} mass of $\mu_n$ to infinity. However, the condition $\mu_{\infty}((0,+\infty))>0$ excludes possibility that {\it all} mass of $\mu_n$ escapes to infinity. By a push-forward under $x\mapsto {1}/(1+x)$, Assumption \ref{ass:A1} is equivalent to weak convergence of the empirical distribution of the random probabilities $\eta_{j:n}=1/(1+\lambda_{j:n})$ to some distribution on $[0,1]$ with positive mass on $(0,1)$.

Another equivalent form of condition \ref{ass:A1} will be using different terminology that is very useful throughout the paper. 
\begin{definition}\label{def:Stieltjes}
The \emph{Cauchy-Stieltjes transforms} of the reflected measure $\mu_n(-\cdot)$ on $[-\infty, 0]$ is defined by
	\begin{equation}\label{eq:phi_n_def_and_conv}
		\Sti_n(z):=\int_{[0,+\infty]} \frac{\mu_{n}({\rm d}x)}{z+x},\qquad z\in\mathbb{C}\setminus (-\infty,0],n\in\mathbb{N}\cup\{\infty\}.
	\end{equation}
where we continuously extend the Cauchy-kernel at $x=\infty$ to $\frac 1{z+\infty}=0$. If $n<\infty$, then $\Sti_n(z)=\frac{P_n'(z)}{nP_n(z)}$. Occasionally in this section and more excessively during the proof, we will also make use of the notation $\Phi_n(z):=z\Sti_n(z)$.
\end{definition}
For every fixed $z\in\mathbb{C}\setminus (-\infty,0]$, the kernel $x\mapsto 1/(z+x)$ is continuous and bounded on $[0,\infty)$, hence assumption \ref{ass:A1} yields
\begin{equation*}
\Sti_n(z)
~\to~ \int_0^\infty\frac{\mu_{\infty}({\rm d}x)}{z+x}=:\Sti_{\infty}(z)\qquad\text{a.s.},\qquad n\to\infty.
\end{equation*}
Moreover, the continuity theorem for Stieltjes transforms, see Proposition 2.1.2(vi) in~\cite{Pastur+Shcherbina}, implies that the above convergence is locally uniform on $\mathbb{C}\setminus (-\infty,0]$ and, in fact, \ref{ass:A1} is equivalent to

\begin{ass}{$({\bf A}_1^{\prime})$}\label{ass:A1'}
For every compact set $K\subset \mathbb{C}\setminus (-\infty,0]$, assume
	\begin{equation}\label{eq:A1'}
\lim_{n\to\infty}\sup_{z\in K}|\Sti_n(z)-\Sti_{\infty}(z)|=0\qquad\text{a.s.}
\end{equation}
\end{ass}

The assumption $\mu_{\infty}((0,+\infty))>0$ implies $\underline{m}<\overline{m}$, where
\begin{equation}\label{eq:um_om_def}
\underline{m}:=\mu_{\infty}(\{0\})
,\qquad \overline{m}:=1-\mu_{\infty}(\{+\infty\})=\mu_\infty(\R).
\end{equation}
The function $x\mapsto x\Sti_{\infty}(x)$ is strictly increasing on $(0,+\infty)$ and
\begin{equation}\label{eq:mu_infty_range}
\lim_{x\to 0+}x\Sti_{\infty}(x)=\mu_{\infty}(\{0\})=\underline{m},\quad \lim_{x\to +\infty}x\Sti_{\infty}(x)=\mu_{\infty}([0,+\infty))=\overline{m}.
\end{equation}
In particular, $\Phi_{\infty}(x):=x\Sti_\infty(x)$ possesses a strictly increasing continuous inverse $\Phi^{\leftarrow}_{\infty}:(\underline{m}, \overline{m})\to (0,+\infty)$, becoming a diffeomorphism. 

\begin{definition}\label{def:logpot}
Let $n\in\mathbb{N}\cup\{\infty\}$ and $\mu_n$ be a probability distribution on $[0,+\infty]$. We shall also use the \emph{normalized logarithmic potential} of its reflection $\mu_n(-\cdot)$, defined by
\begin{equation}\label{eq:psi_def}
\Psi_{n}(z):=\int_0^\infty\log \left(\frac{z+x}{1+x}\right)\mu_{n}({\rm d}x) ,\qquad z\in \mathbb{C}\setminus (-\infty,0] .
\end{equation}
where $\log$ is the principal branch of the complex logarithm, extended to $\log \left(\frac{z+\infty}{1+\infty}\right)=0$. If $n<\infty$, then $\Psi_{n}(z)=\frac 1 n \log P_n(z)$. If the logarithm is $\mu_{n}$-integrable at infinity, enforcing $\mu_n(\{+\infty\})=0$, then we may also define the (classical, non-normalized)  \emph{logarithmic potential} of  $\mu_n(-\cdot)$ by
\begin{equation}\label{eq:logpot_def}
	\U_{n}(z):=\int_0^\infty\log \left( {z+x} \right)\mu_{n}({\rm d}x) ,\qquad z\in \mathbb{C}\setminus (-\infty,0] .
\end{equation} 
If $n<\infty$, then $\U_{n}(z) =\frac 1 n \log p_n(z)$.
\end{definition}
These measure transforms are related as follows
\begin{equation}\label{eq:Phi_Psi_integral}
\Psi_n(z)=\int_1^z {\Sti_n(u)} {\rm d}u=\U_n(z)-\U_n(1), \qquad \Psi_n'(z)=\U_n'(z)=G_n(z),
\end{equation}
for $z\in\mathbb{C}\setminus (-\infty,0]$ and where the integration is taken along the segment\footnote{Instead of the segment $[1, z]$ we may take any sufficiently smooth curve in $\mathbb{C}\setminus (-\infty,0]$ connecting the points $1$ and $z$.} $[1, z]$. 
It will also be important in what follows that the random functions $\Sti_n,\Phi_n,\Psi_n$ and $\U_n$ are a.s.~analytic on $\mathbb{C}\setminus (-\infty,0]$ and, particularly, on the right open half-plane $$\mathbb{H}:=\{z\in\mathbb{C}:\Re(z)>0\}\subset \C.$$

In particular, locally uniform convergence of analytic functions \ref{ass:A1'} is equivalent to
\begin{equation}\label{eq:A1''}
	\lim_{n\to\infty}\sup_{z\in K}|\Psi_n(z)-\Psi_{\infty}(z)|=0\qquad\text{a.s.} 
\end{equation}
for every compact set $K\subset \mathbb{C}\setminus (-\infty,0]$.

The assumption \ref{ass:A1} (or \ref{ass:A1'} or \eqref{eq:A1''}) alone yields the first order asymptotic behavior of $g_n$. The following proposition is an immediate consequence of \cite[Theorem 2.4, 5.1]{JalKabMar:I} in conjunction with the representation~\eqref{eq:coeff}.

\begin{proposition}\label{prop:JalKabMar2024}
Suppose \ref{ass:A1} and let $\underline{m}$, $\overline{m}$ be as defined in~\eqref{eq:um_om_def}. Then, there exists an infinitely differentiable strictly concave function $g:(\underline{m},\overline{m})\to (-\infty,0]$ such that,
for every compact set $K\subset (\underline{m},\overline{m})$,
\begin{equation}\label{eq:uniform_convergence_of_profiles_new}
\sup_{\alpha\in K}\left|g_n(\alpha)-g(\alpha)\right|\ton 0\qquad\text{a.s.},
\end{equation}
where $g_n$ is defined by~\eqref{eq:g_n_def}. The function $-g$ is the Legendre transform of a strictly convex function $u\mapsto \Psi_{\infty}(\eee^u)$, $u\in\mathbb{R}$, where $\Psi_{\infty}$ is defined in~\eqref{eq:psi_def}. More precisely,
\begin{equation}\label{eq:g_def_as_Legendre}
-g(\alpha)=\sup_{u\in\mathbb{R}}(\alpha u-\Psi_{\infty}(\eee^u))=\alpha\log\Phi_{\infty}^{\leftarrow}(\alpha)-\Psi_{\infty}(\Phi_{\infty}^{\leftarrow}(\alpha)),\quad \alpha\in (\underline{m},\overline{m}),
\end{equation}
where $\Phi_{\infty}^{\leftarrow}$ is the inverse function of $\Phi_{\infty}(x)=x\Sti_\infty(x)$ defined in Definition \ref{def:Stieltjes}. Furthermore, the following identity holds for all $\alpha\in (\underline{m},\overline{m})$,
\begin{equation}\label{eq:phi_Phi_inversion}
\eee^{-g'(\alpha)}= \Phi_{\infty}^{\leftarrow}(\alpha).
\end{equation}
\end{proposition}

In order to quantify fluctuations of $g_n$ around $g$ it is natural to impose assumptions on the second order asymptotic behavior of the random measures $\mu_n$, as $n\to\infty$. To this end, let $\mathcal{A}(\mathbb{H})$ be the space of functions which are analytic in $\mathbb{H}$. We endow the space $\mathcal{A}(\mathbb{H})$ with the topology of locally uniform convergence. It is known that this space is a separable completely metrizable topological space (Polish space). Condition~\ref{ass:A1'} implies that the sequence of random analytic functions $(\Sti_n(z))_{z\in\mathbb{H}}$ converges a.s.~in this space to a deterministic analytic function $(\Sti_{\infty}(z))_{z\in\mathbb{H}}$. This can be regarded as a version of the strong law of large numbers. Thus, assumption \ref{ass:A2'} given below, which is a version of the functional central limit theorem quantifying fluctuations of $\Sti_n$ around $\Sti_{\infty}$, is a rather natural condition on the second order behavior of $\mu_n$. However, it will often turn out that the centering $\Sti_\infty$ does not capture deterministic finite $n$  fluctuations, and will thus be replaced by a deterministic sequence $(\widetilde\Sti_n)\subset\mathcal{A}(\mathbb{H})$ converging to $\Sti_{\infty}$.
We begin with the equivalent assumption in terms of normalized logarithmic potentials $\Psi_n$.

\begin{ass}{$({\bf A}_2)$}\label{ass:A2}
	Suppose that there exist a deterministic sequence $(c_n)_{n\in\mathbb{N}}$ of positive numbers with $ c_n\to\infty$, a deterministic sequence of analytic functions $(\widetilde\Psi_n )\subset \mathcal{A}(\mathbb{H})$, which are real on $\R$ with $\widetilde\Psi_n(1)=0$, and a random element $(\intprocleft(z))_{z\in\mathbb{H}}\in\mathcal{A}(\mathbb{H})$ such that
	\begin{equation}\label{eq:A2}
		c_n\left(\Psi_n(z)-\widetilde\Psi_n(z)\right)_{z\in\mathbb{H}}~\toweak~(\intprocleft(z))_{z\in\mathbb{H}}
	\end{equation}
	on $\mathcal{A}(\mathbb{H})$.
\end{ass} 

Since the topology of locally uniform convergence respects differentiation  and integration, the following differentiated version \ref{ass:A2'} is equivalent to \ref{ass:A2}, with	\begin{align}\label{eq:IF}
	\intprocleft(z)=\int_1^z\intprocleft'(u){\rm d}u,\qquad	 
	\widetilde\Psi_n(z)=\int_1^{z} {\widetilde\Sti_n(u)} {\rm d}u,\quad z\in\mathbb{H}.
\end{align}
\begin{remark}
 The subscript $1$ in $\intprocleft$ indicates that the process is pinned at $1$, that is, $\intprocleft(1)=0$ a.s. We reserve the notation $\intprocleftup$ for a companion unpinned process, which will be introduced in Assumption~\ref{ass:A3}.
\end{remark}

 \begin{ass}{$({\bf A}_2^{\prime})$}\label{ass:A2'}
Suppose that there exist a deterministic sequence $(c_n)_{n\in\mathbb{N}}$ of positive numbers with $ c_n\to\infty$, a deterministic sequence of analytic functions $(\widetilde\Sti_n)\subset \mathcal{A}(\mathbb{H})$, which are real on $\R$, and a random element $(\intprocleft'(z))_{z\in\mathbb{H}}\in\mathcal{A}(\mathbb{H})$ such that 
\begin{equation}\label{eq:A2'}
c_n\left(\Sti_n(z)-\widetilde\Sti_n(z)\right)_{z\in\mathbb{H}}~\toweak~(\intprocleft'(z) )_{z\in\mathbb{H}}
\end{equation}
weakly on $\mathcal{A}(\mathbb{H})$.
\end{ass}

As usual, convergence~\ref{ass:A2'} is equivalent to the convergence of finite-dimensional distributions plus tightness in the space of probability measures on $\mathcal{A}(\mathbb{H})$, where the former can be dealt with using the Cram\'{e}r-Wold device. Indeed, according to Proposition 2.5 in~\cite{Shirai:2012} and a remark after Lemma 2.6 therein, a simple sufficient moment condition for tightness is the following $p$-th moment bound, and a version of the Vitali-Porter theorem for random processes, see Theorem~\ref{thm:Vitali-Porter-random} in the Appendix, shows that, assuming tightness, it suffices to check only real positive $t_1,\ldots,t_m$ and real $\gamma_1,\ldots,\gamma_m$. This demonstrates that the following `real' version implies Assumption \ref{ass:A2'}, which  is typically easier to check.
 
  \begin{ass}{$({\bf A}_2^{\prime\prime})$}\label{ass:A2''}
Suppose that there exist a deterministic sequence $(c_n)_{n\in\mathbb{N}}$ of positive numbers with $ c_n\to\infty$, a deterministic sequence of analytic functions $(\widetilde\Sti_n)\subset \mathcal{A}(\mathbb{H})$, which are real on $\R$, and a random element $(\intprocleft'(z))_{z\in\mathbb{H}}\in\mathcal{A}(\mathbb{H})$ such that
	\begin{equation}\label{eq:A2_fd}
		c_n\sum_{i=1}^{m}\gamma_i \left(\Sti_n(z_i)-\widetilde\Sti_n(z_i)\right)~\toweak~\sum_{i=1}^{m}\gamma_i \intprocleft'(z_i) ,
	\end{equation}
	for every $m\in\mathbb{N}$, $\gamma_1,\ldots,\gamma_m\in\mathbb{R}$ and $z_1,\ldots,z_m>0$.  Suppose further that, for some $p>0$, the function 
\begin{align}\label{eq:tightness}
z\mapsto \sup_{n\in\mathbb{N}}\big(c_n^p\mathbb{E}\left[|\Sti_n(z)-\widetilde\Sti_n(z)|^p\right]\big).
\end{align}
is locally integrable on $\mathbb{H}$.
\end{ass}

\smallskip
Any of the above assumptions \ref{ass:A2} together with Assumption \ref{ass:A1}, imply
\begin{equation}\label{eq:centering_to_phi_infty}
	\widetilde\Sti_n(z)\to \Sti_{\infty}(z) ,\quad 	\widetilde\Psi_n(z)\to \Psi_{\infty}(z),\quad \widetilde\Phi_n(z):=z\widetilde{\Sti}_n(z)\to zG_\infty(z)=:\Phi_\infty(z)
\end{equation}
locally uniformly on $\mathbb{H}$ as $n\to\infty$. Moreover, by analyticity $\widetilde\Sti_n\in\mathcal A(\mathbb H)$ we have $\widetilde\Sti_n'\to \Sti_\infty'$ locally uniformly. On any compact set $K^{\prime}\subset (0,+\infty)$ we have $\Phi_\infty'>0$. By the locally uniform convergence this entails $\widetilde\Phi_n'>0$ on $K^{\prime}$ for all $n\geq n_0$, where $n_0=n_0(K^{\prime})$ is sufficiently large. Hence the inverse $\widetilde\Phi_n^\leftarrow$ is well-defined and strictly increasing on any compact $K\subset (\underline{m},\overline{m})$ for all $n\geq n_1$, where $n_1=n_1(K)$ is sufficiently large.  Hence, on any compact $K\subset (\underline{m},\overline{m})$, we may define
\begin{equation}\label{eq:tildeg_n_def_as_Legendre}
-\widetilde g_n(\alpha):=\alpha\log \widetilde\Phi_n^{\leftarrow}(\alpha)-\widetilde\Psi_n(\widetilde\Phi_n^{\leftarrow}(\alpha)),\quad \alpha\in K\subset (\underline{m},\overline{m})
\end{equation} 
for $n$ sufficiently large.

In typical examples, we may choose $\widetilde\Sti_n=\E[\Sti_n]\to\Sti_\infty$, or  even $\widetilde\Sti_n\equiv\Sti_{\infty}$ independent on $n$. In the latter case, $\widetilde\Psi_n\equiv\Psi_{\infty}$, whence $\widetilde g_n\equiv g$ from Proposition \ref{prop:JalKabMar2024}. In the general case, we always have
$$
\widetilde g_n(\alpha)~\to g(\alpha),\text{ locally uniformly on }(\underline{m},\overline{m}),\quad n\to\infty.
$$

\begin{remark}[Legendre representation]\label{rem:not-Legendre}
	Note that the function $u\mapsto \widetilde\Psi_n(\eee^u)$ is strictly convex on compacts (for sufficiently large $n$) and thus, \eqref{eq:tildeg_n_def_as_Legendre} may be viewed as a restricted Legendre transform of $\widetilde\Psi_n(\eee^u)$, analogous to  \eqref{eq:g_def_as_Legendre}. However,  $\widetilde\Phi_n$ may fail to be strictly monotone close to the origin at finite $n$, such that $-\widetilde g_n$ need not coincide with the full Legendre transform $\sup_{u\in\R}(\alpha u-\widetilde\Psi_n(\eee^u))$. The origin does not cause any trouble in the following, since we are interested in the topology of locally uniform convergence.
\end{remark}

\begin{remark}[Notation summary]\label{rem:notation}
We denote the analytic random measure transforms by capital letters $\Sti,\U,\Phi,\Psi\in\mathcal A(\mathbb H)$, as in \cite{JalKabMar:I}, where $\Sti,\U$ are standard notations in the literature, and their modifications $\Phi,\Psi$ will turn out to be most convenient for our proofs and inverses. Recall that the exponential profile $g(\alpha)$ and its approximation $g_n$ defined in \eqref{eq:g_n_def}, denote functions on  $(\underline{m},\overline{m})$. All their deterministic centerings will be equipped with $\widetilde\cdot$, that is $\widetilde g_n,\widetilde\Sti_n,\widetilde\Phi_n,\widetilde\Psi_n,\widetilde\U_n$. Our random fluctuation limits have fraktur letters, such as $\intprocleft(z),\intprocleft'(z)$, or $\Gau(\alpha)$,  with the latter reserved for Gaussian processes in the variable $\alpha$. 
\end{remark}

\begin{remark}\label{rem:cov_I}
	If $\intprocleft'$ is a Gaussian process with covariance function $\kappa_{\intprocleft'}(z,w)$, then $\intprocleft$ is also a Gaussian process, since (pinned) integration $\intprocleft'\mapsto \int_1^{(\cdot)}\intprocleft'(z){\rm d}z=\intprocleft$ is a linear functional. It inherits the covariance function by linearity
	$$ \kappa_{\intprocleft}(z,w)=\int_1^z\int_1^w  {\kappa_{\intprocleft'}(u,v)} \dd u\dd \bar v.$$
\end{remark}
\subsection{The main result for coefficients}
In what follows we denote by $\mathcal C(A,\mathbb{R}^d)$ the space of continuous functions defined on $A$ and taking values in $\mathbb{R}^d$, and by $\mathcal D(A,\mathbb{R}^d)$ the Skorokhod space of c\`{a}dl\`{a}g functions. Both spaces $\mathcal C(A,\mathbb{R}^d)$ and $\mathcal D(A,\mathbb{R}^d)$ are endowed with the topology of locally uniform convergence. If $d=1$, we shall write $\mathcal C(A)$ instead of $\mathcal C(A,\mathbb{R})$ and $\mathcal D(A)$ instead of $\mathcal D(A,\mathbb{R})$. Observe that we equip $\mathcal{D}(A)$ with the topology of locally uniform convergence rather than the usual Skorokhod $J_1$ topology. Since all limiting processes considered below have continuous sample paths, convergence in the Skorokhod $J_1$ topology is equivalent to locally uniform convergence in this setting.

We are ready to formulate our main result.

\begin{theorem}\label{thm:main}
Assume \ref{ass:A1} and \ref{ass:A2} with $c_n =\mathcal O(n)$  as $n\to\infty$ and recall $g$ from Proposition~\ref{prop:JalKabMar2024}, $g_n$ from \eqref{eq:g_n_def} and $\widetilde{g}_n$ from \eqref{eq:tildeg_n_def_as_Legendre}. Then,
\begin{align}\label{eq:main_flt_measures}
	c_n\left(g_n(\alpha)-\widetilde g_n(\lfloor \alpha n \rfloor/n)+\tfrac 1 n  \log \sqrt{\frac{2\pi n}{-g''(\alpha)}}  \right)_{\alpha\in (\underline{m},\overline{m})}~\toweak~\left(\intprocleft(\eee^{-g'(\alpha)})\right)_{\alpha\in (\underline{m},\overline{m})}
\end{align}
on the space of probability measures on $\mathcal D((\underline{m},\overline{m}))$.
\end{theorem}
Strictly speaking, the left hand side of \eqref{eq:main_flt_measures} may be ill-defined simultaneously for all $\alpha \in (\underline m,\overline m)$ or small $n$, since definition \eqref{eq:tildeg_n_def_as_Legendre} of $\widetilde g_n$ requires the threshold for $n\in\N$ to depend on $K\subset (\underline m ,\overline m)$. Similarly, the weak convergence of \eqref{eq:main_flt_measures} should be understood by the very choice of the topology; Choose $K\subset (\underline m,\overline m)$ first, then put $n$ large enough for $\widetilde g_n$ to be defined and send $n\to\infty$ for convergence. Alternatively, set the left hand side of \eqref{eq:main_flt_measures} equal to zero, as long as it is undefined. The proof of Theorem~\ref{thm:main} will be given in Section~\ref{sec:proof}.

\begin{remark}[Alternative presentation]\label{rem:non-linear_time_change} 
	Assuming for simplicity that $\widetilde\Phi_n\equiv\Phi_{\infty}$, one may rewrite \eqref{eq:main_flt_measures} in terms of $\Phi_{\infty}$. Recalling~\eqref{eq:g_def_as_Legendre}, then differentiating \eqref{eq:phi_Phi_inversion} for $\alpha=\Phi_{\infty}(x)$ in~\eqref{eq:main_flt_measures} and using $(\Phi_{\infty}^{\leftarrow})'(\Phi_{\infty}(x))=1/\Phi_{\infty}^{\prime}(x)$, we see that~\eqref{eq:main_flt_measures} is equivalent to
$$
c_n\left(g_n(\Phi_{\infty}(x))+\frac{\lfloor n\Phi_{\infty}(x)\rfloor}{n}\log x - \Psi_{\infty}(x)+\tfrac 1 n\log \sqrt{ {2\pi n x\Phi'(x)}} 
 \right)_{x>0}~\toweak~\left(\intprocleft(x)\right)_{x>0},
$$
on the space of probability measures on $\mathcal D((0,\infty))$. See the proof around \eqref{eq:gtophi} below for details. 
\end{remark}
\begin{remark}[Simplifications]\label{rem:scaling} 
Recall that $g,\widetilde g_n\in\mathcal C^\infty$ with $g''<0$. By the mean value theorem for differentiable functions $g_n(\lfloor \alpha n \rfloor/n)=g_n(\alpha)+O(n^{-1})$. Thus, if the sequence $(c_n)$ in \ref{ass:A2} satisfies $c_n=o(n)$, then $\widetilde g_n(\lfloor \alpha n \rfloor/n)$ can be replaced by $\widetilde g_n(\alpha)$ and the summand $\tfrac 1 n\log \sqrt{\frac{2\pi n}{-g''(\alpha)}}$ in~\eqref{eq:main_flt_measures} can be replaced by $\tfrac 1 n\log\sqrt n$. If moreover $c_n=o(n/\log n)$, then it can be omitted completely.

For triangular arrays $\lambda_{j:n}$ resembling i.i.d.~sequences we will have $c_n\sim\sqrt{n}$, see Section~\ref{sec:independent} below. However, in models coming from random matrix theory one typically has $c_n\sim n$, see Section \ref{sec:randommatrices}.

Lastly, we may replace $\widetilde g_n(\lfloor \alpha n\rfloor /n)$ in \eqref{eq:main_flt_measures} directly by $g(\alpha)$ if we additionally assume $\widetilde\Psi_n=\Psi_\infty+o(c_n^{-1}) $ locally uniformly and $c_n=o(n)$.
Indeed, by our convergences \eqref{eq:centering_to_phi_infty}, both $\widetilde\Phi_n^{\leftarrow}(K)$ and
$\Phi_\infty^{\leftarrow}(K)$ lie in a common compact interval $I\subset(0,\infty)$ for all sufficiently large $n$.
By \eqref{eq:g_def_as_Legendre} and Remark \ref{rem:not-Legendre}, both variational formulas may therefore be restricted to $x\in I$, yielding
$$
c_n\sup_{\alpha\in K}|\widetilde g_n(\alpha)-g(\alpha)|
\le c_n\sup_{x\in I}|\widetilde\Psi_n(x)-\Psi_\infty(x)|
\longrightarrow0.
$$
\end{remark}

\begin{remark}[Heuristic argument of the proof] 
	Let us present a non-rigorous argument, which extends \cite[Remark 2.6]{JalKabMar:I} to the fluctuation scale and which provides an idea of the proof. The log-factor in \eqref{eq:main_flt_measures} appearing for $c_n\sim n$ is quite reminiscent of the Laplace's method, and indeed can formally be explained as follows.
	
	By \eqref{eq:g_n_def}, a Riemann approximation and Laplace's method, we find
\begin{align*}
	P_n(x)= 
	\sum_{k=0}^n \eee^{ng_n(k/n)+k\log x}&\approx n\int_{[0,1]}\eee^{n(g_n(\alpha )+\alpha \log x)}\dd \alpha\\&  \approx   \sqrt{\frac{2\pi n}{- g_n''(\alpha ^*(x))}}\eee^{n(g_n(\alpha ^*(x))+\alpha ^*(x)\log x)},  
	\end{align*}
	where $\alpha ^*(x)$ is a solution to the critical equation $g_n'(\alpha )+\log x=0$. In the large $n$ limit, we expect this equation to be well approximated by $-\widetilde g_n'(\alpha )=\log x$. Since $-\widetilde g_n'=\log(\widetilde\Phi_n^\leftarrow)$ by Legendre inversion of $\widetilde\Psi_n(e^u)$ similar to  \eqref{eq:phi_Phi_inversion}, it follows $\alpha ^*(x)\approx \widetilde\Phi_n(x)$. Plugging in $x=\widetilde\Phi_n^{\leftarrow}(\lfloor \alpha n\rfloor /n )$ for some $\alpha\in (\underline m , \overline m)$ and taking the logarithm, we would get
	\begin{align*} 
			n(g_n(\alpha )-\widetilde g_n(\tfrac{\lfloor \alpha n\rfloor}{n} ))+ \log\sqrt{\frac{2\pi n}{-g''(\alpha )}}
		&\approx\log P_n(\widetilde\Phi_n^{\leftarrow}(\tfrac{\lfloor \alpha n\rfloor}{n} ))-n\big( \widetilde g_n(\tfrac{\lfloor \alpha n\rfloor}{n} )+\tfrac{\lfloor \alpha n\rfloor}{n} \log(\widetilde\Phi_n^{\leftarrow}(\tfrac{\lfloor \alpha n\rfloor}{n} )) \big)\\
		&=\log P_n(\widetilde\Phi_n^{\leftarrow}(\tfrac{\lfloor \alpha n\rfloor}{n} ))-n \widetilde\Psi_n(\widetilde\Phi_n^{\leftarrow}(\tfrac{\lfloor \alpha n\rfloor}{n} ))\\
		&\toweak \intprocleft(\Phi_{\infty}^{\leftarrow}(\alpha ))
	\end{align*}
	using \eqref{eq:tildeg_n_def_as_Legendre} to the second line and the integrated assumption \ref{ass:A2} with $\widetilde\Phi_n^\leftarrow(\lfloor \alpha n\rfloor /n ) \to\Phi_\infty^\leftarrow(\alpha)$ in the last line.
	
	This is the claim of Theorem \ref{thm:main}. Our proof does not explicitly follow the above Laplace heuristic, but instead relies on the probabilistic representation \eqref{eq:P_n}. A local large deviation principle will formalize the above Laplace method and we will see in Lemma \ref{lem:first_summand} that $x=\widetilde\Phi_n^{\leftarrow}(\lfloor \alpha n\rfloor /n )$ becomes the tilting in this LDP, up to a random error.
\end{remark}

\subsection{The main result for symmetric elementary polynomials}

Recall from \eqref{eq:p_n} the monic polynomial $p_n(z):=p_n(1)\cdot P_n(z)$ with roots $-\lambda_{j:n}$ whose decreasing coefficients are the symmetric elementary polynomials $e_{k:n}$ of its (reflected) roots.
We replace the normalized coefficients $a_{k:n}$ of the previous section by $e_{k:n}=p_n(1)a_{n-k:n}$ as in \eqref{eq:coeff} and its profile $g_n(\alpha)$ from \eqref{eq:g_n_def} by 
\begin{align}\label{eq:decomposition_e_to_a}
\tfrac{1}{n}\log e_{\lfloor \alpha n\rfloor:n}=\tfrac{1}{n}\log a_{n-\lfloor \alpha n\rfloor:n}+\tfrac 1 n\log  p_n(1).
\end{align} It is necessary to impose the following stronger condition (compared to Assumption \ref{ass:A2}), which also controls fluctuations of the (unnormalized) logarithmic potential $\U_n(z)=\tfrac 1 n \log p_n(z)$ defined in \eqref{eq:logpot_def}, as well as its normalization $\U_n(1)=\frac 1 n\log p_n(1)$.
\begin{ass}{$({\bf A}_3)$}\label{ass:A3}
	Suppose that there exist a deterministic sequence $(c_n)_{n\in\mathbb{N}}$ of positive numbers with $ c_n\to\infty$, a deterministic sequence of analytic functions $(\widetilde \U_n(z))\subset \mathcal{A}(\mathbb{H})$, which are real on $\R$, and a random element  $(\intprocleftup(z))_{z\in\mathbb{H}}$ of the space $\mathcal{A}(\mathbb{H})$ such that
	\begin{equation}\label{eq:A3}
		c_n\left(\U_n(z)-\widetilde \U_n(z)\right)_{z\in\mathbb{H}}~\toweak~(\intprocleftup (z))_{z\in\mathbb{H}}
	\end{equation}
weakly on $\mathcal{A}(\mathbb{H})$.
\end{ass}

It is easy to see that \ref{ass:A3} implies \ref{ass:A2} for the normalized logarithmic potential $\Psi_n(z)=\U_n(z)-\U_n(1)$ from \eqref{eq:psi_def}, by pinning all the processes at $1$, that is, by setting $\widetilde\Psi_n(z)=\widetilde \U_n(z)-\widetilde \U_n(1)$, $\intprocleft(z):=\intprocleftup(z)-\intprocleftup(1)$ and writing
\begin{equation}\label{eq:Psi_to_L}
c_n\left(\Psi_n(z)-\widetilde \Psi_n(z)\right)=	c_n\left(\U_n(z)-\widetilde \U_n(z)\right)-c_n\left(\U_n(1) -\widetilde \U_n(1)\right)~\toweak~\intprocleft (z).
\end{equation}
In general,~\ref{ass:A2} gives no information of the rightmost normalization term in \eqref{eq:Psi_to_L}, and hence it does not imply \ref{ass:A3}.  Both conditions become equivalent, if the fluctuations vanish at infinity as follows.
 
\begin{proposition}\label{prop:Billing}
The following are equivalent
\begin{itemize} 
	\item Condition \ref{ass:A2}, and additionally, for some $\widetilde \U_n$ satisfying $\widetilde\Psi_n(z)=\widetilde \U_n(z)-\widetilde \U_n(1)$ and every $\varepsilon>0$, we have
	\begin{equation}\label{eq:B2}
		\lim_{x\to+\infty}\limsup_{n\to\infty}\mathbb{P}\left[c_n\left| \U_n(x)-\widetilde \U_n(x) \right|>\varepsilon\right]=0.
	\end{equation}
	\item Condition \ref{ass:A3}, and additionally $ \intprocleftup (x)\to 0$ in probability as $x\to \infty$.
\end{itemize}
\end{proposition} 
Note that if $\U_n(z)$ even remains bounded in probability, then we automatically have $\overline m=1-\mu_\infty(\{+\infty\})=1$. 
Finally, let us present the implication of \ref{ass:A3} for symmetric elementary polynomials. Assuming \ref{ass:A1} and \ref{ass:A3}, the function $\widetilde g_n$ from \eqref{eq:tildeg_n_def_as_Legendre} takes the form
\begin{equation} \label{eq:tildeg_n_def_as_Legendre2}
	-\widetilde g_n(\alpha):=\alpha\log\left((z\widetilde \U_n'(z))^\leftarrow(\alpha)\right)-\widetilde \U_n((z\widetilde \U_n'(z))^\leftarrow(\alpha))+\widetilde \U_n(1),\quad \alpha\in K
\end{equation} 
for any compact $K\subset (\underline{m},\overline{m})$ and $n$ sufficiently large.

\begin{corollary}\label{cor:main_symm}
	Suppose \ref{ass:A1} and \ref{ass:A3} with $c_n=\mathcal O (n)$, then
	\begin{align*}
		c_n\left( \tfrac{1}{n}\log (e_{\lfloor \alpha n \rfloor:n}) -\widetilde g_n\big(1-\tfrac{\lfloor \alpha n\rfloor}{n}\big)-\widetilde \U_n(1)+\tfrac 1 n  \log \sqrt{\frac{2\pi n}{-g''(1-\alpha )}} \right) \toweak~\intprocleft (\eee ^{-g'(1-\alpha)}) 
	\end{align*}
	on the space of probability measures on $\mathcal D((1-\overline m,1-\underline{m}))$.
\end{corollary}
This corollary follows directly from Theorem \ref{thm:main}, together with the decomposition \eqref{eq:decomposition_e_to_a} and the above discussion of \eqref{eq:Psi_to_L}. A short proof of Proposition \ref{prop:Billing} will be given after the proof of Theorem \ref{thm:main}.

\section{Applications}\label{sec:appl}
Before we come to the proof of Theorem \ref{thm:main}, this section will deal with three applications in increasing difficulty: First, completely dependent roots, second, independent roots and third, dependent roots of characteristic polynomials of random matrices.

\subsection{Completely dependent roots and prescribed independent coefficients} 
As a warm-up, we begin with a particularly elementary example of Theorem \ref{thm:main}, which can also be verified directly using the binomial theorem and Stirling’s formula. Besides providing a sanity check, it leads to a construction of real-rooted random polynomials whose logarithmic coefficients have independent Gaussian limiting fluctuations at prescribed locations. This complements classical random polynomial models with independent coefficients, which generically have complex roots.
\begin{example}
For $b>0$, $c_n\to\infty$ with $c_n=o(n/\log n)$ and $Z\sim\mathcal N (0,1)$, define
	\begin{align}
		\lambda_{j:n}:= b\exp\left( \frac {Z}{bc_n}\right),\quad j=1,\dots, n.
	\end{align}
Then, the coefficients $a_{k:n}$ of $P_n(z) =\prod_{j=1}^{n}\frac{z+\lambda_{j:n}}{1+\lambda_{j:n}}=\left(\frac{z+b\eee^{Z/(bc_n)}}{1+b\eee^{Z/(bc_n)}}\right)^n$ satisfy
$$
a_{\lfloor \alpha n\rfloor:n}=(1+b\eee^{Z/(bc_n)})^{-n}\binom{n}{\lfloor \alpha n\rfloor}(b\eee^{Z/(bc_n)})^{n-\lfloor \alpha n\rfloor}
$$
and asymptotically
$$
	c_n\left(\frac{1}{n}\log a_{\lfloor  \alpha n\rfloor:n}- \log\left(
	\frac{b^{ 1-\alpha}}
	{(1+b) \alpha^\alpha(1-\alpha)^{1-\alpha}}
	\right)\right)_{\alpha\in(0,1)}
	~\toweak~
\frac{Z}{b}\left(\frac1{1+b}-\alpha\right)_{\alpha\in(0,1)} .
	$$
\end{example}
Note that the roots satisfy $\lambda_{j:n} =b+\frac Z{c_n}+O_{\mathbb P}(c_n^{-2})$ and $(\lambda_{j:n})_j$  are identically dependent through a `shift' realization $Z/c_n$. Thus, also the coefficients are strongly correlated through $Z$. 
\begin{proof}
	Observe that $\mu_\infty=\delta_b$ with
$$
	\Sti_n(z)=\frac{1}{z+b\exp(Z/(bc_n))}
	\longrightarrow
	\Sti_\infty(z)=\frac{1}{z+b}
$$
	almost surely locally uniformly on $\mathbb H$.
Solving $x\Sti_\infty(x)=\alpha$ gives $\eee^{-g'(\alpha)}=b\alpha/(1-\alpha)$, and
	$\widetilde\Sti_n=\Sti_\infty$ satisfies Assumption \ref{ass:A2'}. Indeed, Taylor expansion gives, almost surely
	locally uniformly on $\mathbb H$,
$$
	\Sti_n(z)
	=\frac{1}{z+b}
	-\frac{Z}{c_n(z+b)^2}
	+\mathcal O(c_n^{-2}),
$$
	and hence
$$
	c_n\bigl(\Sti_n(z)-\widetilde\Sti_n(z)\bigr)
	\longrightarrow -\frac{Z}{(z+b)^2}.
$$
	By integration, we obtain
	$\intprocleft(z)=Z /(z+b)-Z/(1+b)$. Then,
	Theorem \ref{thm:main} and the simplification for $c_n=o(n/\log n)$ discussed in Remark \ref{rem:scaling} imply
	\begin{align}\label{eq:bibedibabedibu}
		c_n\bigl(g_n(\alpha)-\widetilde g_n(\alpha)\bigr)
		\toweak
		\intprocleft(\eee^{-g'(\alpha)})
		=\frac Zb\left(\frac1{1+b}-\alpha\right)
	\end{align}
	in $\mathcal D((0,1))$. With our deterministic centering, we have $\widetilde\Psi_n(z)=\Psi_\infty(z)=\log\frac{z+b}{1+b}$. Thus, the definition \eqref{eq:tildeg_n_def_as_Legendre} gives exactly
	\begin{align*}
		\widetilde g_n(\alpha)
		 =-\alpha\log\alpha
		+(1-\alpha)\log\frac{b}{1-\alpha}
		-\log(1+b)
		=g(\alpha).
	\end{align*}
The claim follows from \eqref{eq:bibedibabedibu}, which also resembles the replacement of $\widetilde g_n$ by $g$ as discussed in the end of Remark \ref{rem:scaling}, only that it is exact here.
\end{proof}
	
One may straightforwardly generalize the above example to independent coefficients at prescribed locations $\alpha_1,\dots,\alpha_m$ for some fixed number $m\in\N$. Throughout this example assume that $n$ runs through $m\mathbb{N}$, the set of positive multiples of $m$. Divide the index set $ \{1,\dots,n\}$ into $m\in\N$ equal-sized blocks $S_k=\{(k-1)\frac{n}{m}+1,\dots, k\frac{n}{m}\}$, $k=1,2,\ldots,m$. We prescribe the distinct limiting roots $b_1,\dots,b_m>0$ with $\mu_\infty=\tfrac 1 m\sum_{k=1}^m\delta_{b_k}$ and
$$
\Sti_\infty(z)=\frac1m\sum_{k=1}^m\frac{1}{z+b_k},
$$ 
but we will couple them by some dependent Gaussian $Z_1 ,\dots,Z_m$. Fix distinct $\alpha_1,\dots,\alpha_m\in(0,1)\setminus\{\Sti_\infty(1)\}$ and define the $m\times m $-matrix $B$ by its entries
$$
B_{jk}=\frac1m\left(\frac1{\eee^{-g'(\alpha_j)}+b_k}-\frac1{1+b_k}\right)=\frac1m\left(\frac{1-\eee^{-g'(\alpha_j)}}{(\eee^{-g'(\alpha_j)}+b_k)(1+b_k)}\right),
$$ 
where $\alpha\mapsto \eee^{-g'(\alpha)}$ is the inverse of $z\mapsto \frac{z}{m}\sum_{k=1}^m\frac{1}{z+b_k}$. The matrix $B$ is invertible by the Cauchy determinant formula, since the $b_k$ are pairwise distinct, the numbers $e^{-g'(\alpha_j)}$ are pairwise distinct by strict monotonicity of $\Phi_{\infty}(x)=xG_\infty(x)$, and $\alpha_j\ne G_\infty(1)$ ensures that $\eee^{-g'(\alpha_j)}\ne1$.

Set $Z \sim\mathcal   N_m(0,(B^{\mathsf T}B)^{-1})$  and choose roots
$$\lambda_{j:n}
	=b_k\exp\left(\frac{Z_k}{b_kc_n}\right),
	\qquad j\in S_k,\quad 1\leq k\leq m,
$$
Applying the calculation from the preceding example to each
block and using $|S_k|/n=1/m$ verifies
Assumption \ref{ass:A2'} with $\widetilde\Sti_n=\Sti_\infty$,
$c_n=o(n/\log n)$, $\widetilde\Psi_n=\Psi_\infty$,
$\widetilde g_n=g$ and integrated limiting process
$$
\intprocleft(z)
=\frac1m\sum_{k=1}^m Z_k
\left(\frac1{z+b_k}-\frac1{1+b_k}\right).
$$ Theorem \ref{thm:main}, together with
Remark \ref{rem:scaling}, therefore showed 
\begin{align}
c_n\left(
\frac1n\log a_{\lfloor\alpha n\rfloor:n}-g(\alpha)
\right)_{\alpha\in(0,1)}
\toweak \Gau(\alpha):=\intprocleft(e^{-g^{\prime}(\alpha)}).\label{eq:indepcoeff1}
\end{align}
as claimed in Application 3 from the introduction.
In particular,
\begin{align}
	c_n\left(
	\frac1n\log a_{\lfloor\alpha_j n\rfloor:n}-g(\alpha_j)
	\right)_{j=1}^m
	\toweak \left(
	\intprocleft\bigl(\eee^{-g'(\alpha_j)}\bigr)
	\right)_{j=1}^m
	=BZ \sim\mathcal N_m(0,I_m).\label{eq:indepcoeff}
\end{align}
Hence the limiting fluctuations at the coefficient locations
$\alpha_1,\dots,\alpha_m$ are independent standard Gaussian variables. Observe that we can first choose any collection of pairwise distinct $\alpha_1,\ldots,\alpha_m\in(0,1)$ and then adjust $b_1,\ldots,b_m$ such that $G_{\infty}(1)\neq \alpha_j$, for all $j=1,\ldots,m$. For a simulation of the finite-$n$ process with three locations $\{\alpha_1,\alpha_2,\alpha_3\}=\{1/4,2/4,3/4\}$ and with blocks $b=(1/100,2,100)$, see Figure \ref{fig:sim}.

\begin{remark}
The above model constructs real rooted polynomials whose logarithmic coefficients become independent at prescribed locations. On the other hand, analyticity of the process $\intprocleft$ prevents for instance white noise limits corresponding to independent coefficients $a_{\lfloor \alpha n\rfloor:n}$. Classical models of random polynomials with independent coefficients generically have complex roots, such as Kac \cite{Kac} and Littlewood–-Offord \cite{LO} polynomials or more general models studied, for example, in~\cite{HHJK,KKM,KZ14}.
\end{remark}

\subsection{Independent roots}\label{sec:independent}
We are interested in an application of Corollary \ref{cor:main_symm} to a triangular array $(\lambda_{j:n})_{j=1,\dots,n}$ of independent random variables on $[0,\infty)$. The primary example we have in mind are i.i.d. roots, but also i.i.d. perturbations of deterministic roots are covered in the following.

We assume a logarithmic moment condition
\begin{align}\label{eq:log-moments}
\sup_{n\in\N}\frac 1 n \sum_{j=1}^n\E [(\log(1+\lambda_{j:n}))^{2+\delta}]<\infty
\end{align}
for some $\delta>0$. We further assume that 
\begin{align}\label{eq:Emu_n}
	\E \mu_n \toweak \mu_\infty
\end{align}
for some deterministic $\mu_\infty$ on $[0,\infty)$ with $\mu_\infty(0)<1$ and, for some kernel $\kappa$ and any $z,w>0$, let
\begin{align}\label{eq:cov-conv}
\frac 1 n \sum_{j=1}^n\Cov\big(\log(z+\lambda_{j:n}),\log(w+\lambda_{j:n})\big)\ton \kappa(z,w).
	\end{align}
Note that \eqref{eq:log-moments} and \eqref{eq:Emu_n}  imply
$$ \frac 1 n \sum_{j=1}^n\E\big[\log(z+\lambda_{j:n})\log(w+\lambda_{j:n})\big]\ton \int\log(z+x)\log(w+x) \dd\mu_\infty(x),$$
such that \eqref{eq:cov-conv} is equivalent to the convergence of $\frac 1 n \sum_j \E\log(z+\lambda_{j:n})\E\log(w+\lambda_{j:n})$. 
Let us mention two well studied special cases of the point process of $\lambda_{j:n}$.
\begin{example}\label{examples}\begin{enumerate}
		\item[(i)] Let $r_{j:n}$ be deterministic non-negative points (think of roots of a sequence of deterministic polynomials) whose empirical distribution converges
		$$\rho_n=\frac 1 n \sum_{j=1}^n \delta_{r_{j:n}}\toweak \rho_\infty.$$ Consider the perturbation $\lambda_{j:n}:=r_{j:n}+\eps_{j:n}$, where \ $\eps_{j:n}\geq 0$ are i.i.d.~copies of $\eps\sim \nu$. Then, by a version of the Glivenko-Cantelli theorem for triangular arrays, see~\cite{Wellner}, we have 
		$$
		\frac{1}{n}\sum_{j=1}^{n}\delta_{(r_{j:n},\eps_{j:n})}~\toweak~\rho_{\infty}\otimes \nu
		$$
		and thereupon
		$$
		\frac{1}{n}\sum_{j=1}^{n}\delta_{r_{j:n}+\eps_{j:n}}~\toweak~\rho_{\infty}\ast \nu=:\mu_{\infty}.
		$$
		The condition~\eqref{eq:log-moments} is satisfied if $$\E[(\log(1+\eps))^{2+\delta}]<\infty,\quad \sup_{n\in\N}\int (\log(1+r))^{2+\delta}\dd\rho_n(r)<\infty.$$
		 This also implies \eqref{eq:cov-conv} with 
		\begin{align}\label{eq:cov_iid-perturb}\kappa(z,w)=\int \Cov\big(\log(z+r+\eps),\log(w+r+\eps)\big)\dd \rho_\infty(r)\end{align}
		\item[(ii)] Simplifying further to $r_{j:n}\equiv 0$, we end up with the setting of i.i.d.~copies  $\lambda_{j:n}$ of $\lambda\sim \nu=\mu_\infty$ with $\mu_\infty(\{0\})<1$, finite logarithmic $(2+\delta)$-moment, and 
		\begin{align}\label{eq:cov_iid}\kappa(z,w)=\Cov\big(\log(z+\lambda),\log(w+\lambda)\big).\end{align}
	\end{enumerate}
\end{example}

\begin{theorem}\label{thm:independent}
		Assume \eqref{eq:log-moments}, \eqref{eq:Emu_n} and \eqref{eq:cov-conv} for the distribution of independent roots $-\lambda_{j:n}\le 0$ of the monic polynomial $p_n(z)=\prod_{j=1}^n(z+\lambda_{j:n})$ with $e_{k:n}=[z^{n-k}]\left(p_n(z)\right)$.
		Then, with $g$ as in \eqref{eq:g_def_as_Legendre}, and $\widetilde \U_n(z)=\frac 1 n \sum_{j=1}^n \E\log(z+\lambda_{j:n})$ and $\widetilde{g}_n$ from \eqref{eq:tildeg_n_def_as_Legendre2},
		\begin{align}\label{eq:independent}
		\sqrt n	\left( \tfrac 1 {n} \log (e_{\lfloor \alpha n \rfloor:n}) -\widetilde g_n\big(1-\alpha\big)-\widetilde \U_n(1)\right) \toweak~\Gau(\alpha)
		\end{align} 
		on the space of probability measures on $\mathcal D((0,1-\underline m))$, where $\Gau(\alpha)= \intprocleftup (\eee ^{-g'(1-\alpha)})$ is a centered Gaussian process with covariance 
		\begin{align*}
			{\rm Cov}[\Gau(\alpha),\Gau(\beta)]=\kappa(\eee ^{-g'(1-\alpha)},\eee ^{-g'(1-\beta)}),\quad \alpha,\beta\in(0,1-\underline m).
		\end{align*}
\end{theorem}

\begin{proof}
	Let us first check Assumption \ref{ass:A1}.  
	For any fixed bounded continuous test-function $f\in\mathcal C_b([0,\infty))$, Hoeffding's inequality and the Borel-Cantelli lemma imply that, almost surely,
	$$\int f\dd\mu_n- \int f\dd\E \mu_n=\frac 1 n\sum_{j=1}^n\big(f(\lambda_{j:n})-\E f(\lambda_{j:n})\big) \to 0.
	$$
	Together with \eqref{eq:Emu_n}, this proves $\int f\dd \mu_n\to\int f \dd\mu_\infty$ almost surely, for a fixed $f\in\mathcal C_b([0,\infty))$. This observation alone is not yet sufficient for almost sure weak convergence as required in \ref{ass:A1} because exceptional sets might be different for different $f\in \mathcal C_b([0,\infty))$. However, choosing countably many $f$ that determine weak convergence (such as piecewise linear functions at rational points), we infer the almost sure weak convergence in \ref{ass:A1}. Note that we have $\underline m=\mu_\infty(\{0\})<1$ and $\overline{m}=1$.
	
	We shall now check Assumption \ref{ass:A3} for $c_n=\sqrt n$ and
	$
	\widetilde \U_n(z)=\frac 1 n \sum_{j=1}^n \E\log(z+\lambda_{j:n}).
	$
	Assumption \eqref{eq:log-moments} implies the Lindeberg condition for the triangular array of random variables 	
	\begin{align*}
		X_{j,n}(z):= \log(z+\lambda_{j:n})-\E\log(z+\lambda_{j:n}),\quad z>0.
	\end{align*} 
	Together with the asymptotic average variance from condition \eqref{eq:cov-conv} for $z=w$, all conditions for Lindeberg's Central Limit Theorem for $X_{j,n}(z)$ are fulfilled, see \cite[Theorem 27.2]{Billingsley}, giving
	\begin{align}\label{eq:Lindeberg_CLT}
		\frac 1 {\sqrt n}\sum_{j=1}^nX_{j,n}(z)\toweak \intprocleftup (z)\sim\mathcal N(0,\kappa(z,z)).
	\end{align}
	In accordance with the Cram\'{e}r-Wold device, precisely as in \eqref{eq:A2_fd}, we fix $m\in\mathbb{N}$, $\gamma_1,\ldots,\gamma_m\in\mathbb{R}$ and $z_1,\ldots,z_m>0$ and we get similarly to \eqref{eq:Lindeberg_CLT},
	\begin{align}\label{eq:Cramer-Wold-log}
		\sqrt n &\sum_{i=1}^m\gamma_i\left( \U_n(z_i)-  \widetilde \U_n(z_i)\right)
		= \sum_{i=1}^m \frac {\gamma_i} {\sqrt n}\sum_{j=1}^nX_{j,n}(z_i)\nonumber\\
		&\toweak \sum_{i=1}^m\gamma_i \intprocleftup (z_i)\sim\mathcal N\Big(0,\sum_{i,k}\gamma_i\gamma_k\kappa(z_i,z_k)\Big).
	\end{align}
	Analogously to \eqref{eq:tightness}, tightness of \eqref{eq:Lindeberg_CLT} on $\mathcal A(\mathbb H)$ follows from \eqref{eq:log-moments} yielding a locally uniform bound on 
	$$\sup_{n\in\N}\E\Big\lvert\sqrt n \left( \U_n(z )-  \widetilde \U_n(z)\right)\Big\rvert^2
	=\sup_{n\in\N}\frac 1 n \sum_{j=1}^n\E\big\lvert X_{j,n}(z)\big\rvert^2<\infty.$$
	In conclusion, we have verified Assumption \ref{ass:A3} with $\intprocleftup$ being a centered Gaussian process with covariance $\kappa(z,w)$. Therefore, Corollary \ref{cor:main_symm} applies, which simplifies to the claim according to Remark \ref{rem:scaling}.
\end{proof}

For the Examples \ref{examples}, the statement of Theorem \ref{thm:independent} can be progressively simplified, and, for a change, let us switch to the representation from Remark \ref{rem:non-linear_time_change}.
\begin{corollary}\label{cor:iid-perturbations}
In the setting of Example \ref{examples} (i), where $\lambda_{j:n}$ are i.i.d.-perturbed roots with deterministic limiting distribution $\rho_\infty$ and perturbation distribution $\nu$, we have $$\widetilde \Sti_n(z)=\int \frac{1}{z+x}(\rho_n\ast\nu)(\dd x)\to \Sti_\infty(z)=\int \frac{1}{z+x}(\rho_\infty\ast\nu)(\dd x)$$ locally uniformly, and
\begin{align*}
\sqrt n	\left( \tfrac 1 n \log (e_{n-\lfloor nx\widetilde \Sti_n(x) \rfloor:n}) +x\widetilde \Sti_n(x)\log x- \widetilde \U_n(x) \right) \toweak~ \intprocleftup (x)
\end{align*} on the space of probability measures on $\mathcal D((0,\infty))$, where $\intprocleftup$ is a centered Gaussian process with covariance given by \eqref{eq:cov_iid-perturb}.
\end{corollary}
\begin{corollary}\label{cor:iid}
Consider i.i.d.~copies $\lambda_{j:n}$ of $\lambda\sim\mu_\infty$ as in Example \ref{examples} (ii). Then, $\widetilde \U_n(z)=\E\log(z+\lambda)$ and $\Sti_\infty(z)=\E[\tfrac {1}{z+\lambda}]$ are independent of $n$ and, consequently
	\begin{align*}
		\sqrt n \left( \tfrac 1 n \log (e_{n-\lfloor n \E[\frac {x}{x+\lambda}] \rfloor:n}) + \E\Big[\frac {x\log x}{x+\lambda}\Big]- \E\big[\log(x+\lambda)\big] \right) &\toweak~\intprocleftup (x), 
	\end{align*} on the space of probability measures on $\mathcal D((0,\infty))$.
Equivalently, using in the original representation \eqref{eq:main_flt_measures} in the coordinate $\alpha$, and noting that \eqref{eq:g_def_as_Legendre} implies $$(1-\alpha)g'(1-\alpha)+\E[\log(\eee^{-g'(1-\alpha)}+\lambda)]=g(1-\alpha)+\E[\log(1+\lambda)],$$
	 it follows 
		\begin{align*}
		\sqrt n	\left( \tfrac 1 n \log (e_{\lfloor \alpha n \rfloor:n}) - g(1-\alpha)-\E[\log(1+\lambda)]\right) &\toweak~\Gau(\alpha)=\intprocleftup (\eee^{-g'(1-\alpha)})
	\end{align*}
 weakly on $\mathcal D((0,1-\underline m))$. Here, $\intprocleftup$ is a centered Gaussian process with covariance kernel given by \eqref{eq:cov_iid} and $\Gau$ has covariance
	\begin{align*}\Cov(\Gau(\alpha),\Gau(\beta))=\Cov\big(\log(\eee^{-g'(1-\alpha)}+\lambda),\log(\eee^{-g'(1-\beta)}+\lambda)\big).\end{align*} 
\end{corollary}
See Figure \ref{fig:sim} for finite-$n$ realizations of the process $\Gau$ for $\lambda=U^2$ with $ U\sim\mathrm{Unif}(0,1)$.

By checking the covariance kernel, it can be shown that the process $\intprocleftup$ coincides with the stochastic integral 
$$
\intprocleftup (x)=\int_0^1\log (x+Q_{\infty}^{\leftarrow}(t)){\rm d}W(t),\quad x>0,
$$
where $Q_{\infty}^{\leftarrow}$ is the quantile function of $\lambda$ and $W$ is the standard Brownian bridge. Simplifying further, we end up with fluctuations of binomial coefficients.

\begin{example}[of Corollary \ref{cor:iid} on Example \ref{examples} of Corollary \ref{cor:main_symm} of Theorem \ref{thm:main}]\label{ex:Bernoulli}
	Let $\lambda_{j:n}\sim \mathrm{Ber}(p)$ be i.i.d. and $N_n=\sum_{j=1}^n\lambda_{j:n}\sim\mathrm{Bin}(n,p)$, then $p_n(z)=z^{n-N_n}(1+z)^{N_n}$ and $e_{k:n}=\binom{N_n}{k}$. Corollary \ref{cor:iid} then provides fluctuations for $\frac 1 {\sqrt n}\log \binom{N_n}{\lfloor\alpha n\rfloor}$, whose pointwise analogue follows again directly from Stirling's formula. In our setting, we have $\underline m=1-p$,  $x\Sti_\infty(x)=1-\frac{p}{x+1}$, $\eee^{-g'(1-\alpha)}=\frac{p-\alpha}{\alpha}$ and Corollary \ref{cor:iid} yields
	$$\sqrt n \left( 
	\frac 1 n \log\binom{N_n}{\lfloor\alpha n\rfloor}
	- \bigl[p\log p-\alpha\log\alpha
	-(p-\alpha)\log(p-\alpha)\bigr] 
	\right) 
	\toweak
	\sqrt{p(1-p)}Z\log\frac{p}{p-\alpha} 
	$$
	in $\mathcal D((0,p))$, where $Z\sim\mathcal N(0,1)$.
\end{example}

\begin{remark}[Related Literature]\label{rem:literature_iid}
Let us compare Theorem~\ref{thm:independent} and Corollary \ref{cor:iid} to previous results in the literature.	Analogous to our setting of normalized elementary symmetric polynomials of degree $k_n\sim\alpha n$ of i.i.d. $\lambda_{j:n}$, a scalar CLT was proved in~\cite{Szekely:1982}, assuming strictly positive $\lambda>0$ and $\mathbb{P}[\lambda=0]=0$, and extended in~\cite{vanEs:1986} to non-negative $\lambda$ (such as the Bernoulli-Example \ref{ex:Bernoulli}).  For $k_n\sim\alpha n$ and signed i.i.d. roots $\lambda_{j:n}\in\R$, a scalar CLT analogous to Corollary \ref{cor:iid} was proven in \cite{major_symm_iid} under the assumptions $\mathbb E|\lambda|<\infty$, $\mathbb E[|\lambda|^{-1}\mathbf1_{\{\lambda \ne0\}}]<\infty$ and more. For different regimes, see also 
\cite{van1988elementary} for $k_n=\mathcal O(\sqrt n)$ and \cite{rempala2005approximation} for fixed $k$ and allowing non-identical distributions. Further, a recent contribution \cite{Bochi+Iommi+Ponce:2022} studies stationary ergodic sequences $\lambda_{j:n}$ (not necessarily independent) and proves an almost sure LLN of $\bigl(e_{k_n:n}/\binom{n}{k_n}\bigr)^{1/k_n}$ with limit being the Hal\'{a}sz-Sz\'{e}kely barycenter.

For our setting of non-negative inputs, assumption \eqref{eq:log-moments} substantially weakens the integrability requirements of  \cite{major_symm_iid}. Further, our results strengthen the above scalar CLT's to functional limit theorems of the process, identify its covariance explicitly, and generalize to independent, non-identically distributed triangular arrays.
\end{remark}

\subsection{Eigenvalues of random matrices}\label{sec:randommatrices} 
In the application of this section we will study characteristic polynomials of random matrices, whose roots are eigenvalues with intensely studied repulsive behavior. We will derive surprisingly simple Gaussian fluctuations for the coefficients of the characteristic polynomial. Assumption \ref{ass:A2} follows from CLT's of linear smooth statistics of the eigenvalues. For the ease of presentation, we stick to the classical model of sample covariance matrices built from iid entries with  vanishing fourth-cumulant $\mathbb{E}|X^4|-|\mathbb{E}X^2|^2-2(\mathbb{E}|X|^2)^2=0$. CLT's have been obtained for similar and increasingly general frameworks, such as real matrices, matrices with variance profile, entries with non-zero fourth cumulant\footnote{at the expense of a more involved limiting variance}, invariant ensembles or determinantal Coulomb gases, see \cite[\S 7]{Pastur+Shcherbina}, \cite[\S 9.7]{BS10spectral}, \cite{Sh11CLT}, \cite{najim2016gaussian}, \cite{LP09CLT} and \cite{Joh98CLT}. 

Let ${\bf X}=[X_{i,j}^{(n)}]_{1\leq i\leq n,1\leq j\leq m}$ be an $n\times m$ random matrix with complex entries $X_{ij}^{(n)}$ being i.i.d.~copies of $X$ satisfying
$$
\mathbb{E}[X ]=0,\quad \mathbb{E}[X ^2]=0,\quad \mathbb{E}[|X |^2]=\sigma^2>0,\quad \E [|X |^4]=2\sigma ^4,
$$
for instance Gaussians. Further assume that $m=m(n)$ is such that 
\begin{equation}\label{eq:c_n_wishart}
	\lim_{n\to\infty}\frac{n}{m(n)}=c\in (0,\infty).
\end{equation}
The random matrix 
\begin{equation}\label{eq:Wishart_def}
	{\bf Y}_n:=\tfrac 1 {m(n)}{\bf X}{\bf X}^{*}
\end{equation} 
is called \emph{sample covariance matrix} (or {\it Wishart matrix}, if $X$ is Gaussian). Let $(\lambda_{j:n})_{1\leq j\leq n}\subset [0,\infty)$ be the family of its eigenvalues and 
$$p_n(-z)={\rm det}({\bf Y}_n-z{\bf I}_n)=\prod_{j=1}^{n}(\lambda_{j:n}-z)=\sum_{k=0}^{n}e_{n-k:n}(-z)^k$$
be its characteristic polynomial. 

It is well known that the condition \ref{ass:A1} holds with $\mu_{\infty}={\sf MP}_{\sigma^2,c}$ 
being the Marchenko-Pastur distribution, see for instance \cite[\S3 \& \S6]{BS10spectral}, which is defined by
\begin{align}
\dd {\sf MP}_{\sigma^2,c} (x)&:=\underline{m}\delta_0(\dd x)+\frac{\sqrt{(c_{+}-x)(x-c_{-})}}{2\pi \sigma ^2 c x}\1[x\in (c_{-},c_{+})]\dd x,\\
\text{for }c_{\pm}&:=\sigma^2(1\pm \sqrt{c})^2,\quad \underline{m}:=\max(0,1-1/c).\label{eq:c+}
\end{align}
It has mass $\underline{m}>0$ at the origin iff $c>1$ and satisfies $\overline{m}=1$.

In the notation of Assumption \ref{ass:A3}, we center $\U_n(z)=\tfrac 1 n \log p_n(z)$ by 
\begin{align}
\widetilde \U_n(z):=\int_0^\infty\log(z+x){\sf MP}_{\sigma^2,n/m}({\rm d}x)
	\ton 
	\int_0^\infty\log(z+x){\sf MP}_{\sigma^2,c}({\rm d}x) \label{eq:MPlogpot}
\end{align}
for $z\in\C\setminus(-\infty,0]$ where also the principal branch of $\log$ is defined. This a consequence of the weak continuity of the Marchenko-Pastur distribution in its second parameter together with the fact that the supports of ${\sf MP}_{\sigma^2,n/m}$ remain uniformly bounded.

The exponential profile $g$ as defined in Proposition~\ref{prop:JalKabMar2024}, was explicitly derived in \cite[Theorem 6.1]{JalKabMarII} satisfying
\begin{align}
	g ^{\prime}(\alpha) &=\log\left(\frac{1-\alpha }{\sigma ^2\alpha(c\alpha +1-c)}\right),\label{eq:g'}\\
-g''(\alpha)&= \frac 1 \alpha +\frac 1 {1-\alpha} +\frac c {c\alpha+1-c},\label{eq:g''} 
\end{align}
for $\alpha\in (0,\underline m)$.
Corollary~\ref{cor:main_symm} yields the following functional limit theorem for the process of  coefficients of the characteristic polynomial of $\mathbf{Y}_n$.

\begin{theorem}\label{thm:FLT_RMT}
	Let ${\bf Y}_n$ be a random covariance matrix as defined above as well as $e _{k:n}:=[z^{n-k}]{p_n(z)}$ the elementary polynomials of its eigenvalues (or equivalently, the decreasing coefficients of the reflected characteristic polynomial $p_n$).  
	Then, with $g,\widetilde \U_n$ as above and $\widetilde{g}_n$ from \eqref{eq:tildeg_n_def_as_Legendre2},
	\begin{align}\label{eq:FLT_RMT} 
 \left(  \log e_{\lfloor \alpha n \rfloor:n} -n\widetilde g_n\big(1-\tfrac{\lfloor \alpha n\rfloor}{n}\big)-n\widetilde \U_n(1)+ \log \sqrt{\frac{2\pi n}{-g''(1-\alpha )}} \right)_{\alpha\in (0,\min(1,1/c))} \toweak~(\Gau(\alpha))_{\alpha\in (0,\min(1,1/c))}
\end{align} 
on the space of probability measures on $\mathcal D((0,\min(1,1/c)))$, where $\Gau(\alpha)=\intprocleftup (\eee ^{-g'(1-\alpha)}) $ is a centered Gaussian process with covariance 
\begin{align*}
	{\rm Cov}[\Gau(\alpha),\Gau(\beta)]=- \log\big( 1-c \alpha\beta \big),\quad \alpha,\beta\in(0,\min(1,1/c)).
\end{align*}
Note that \eqref{eq:g''} makes the last term of the centering in \eqref{eq:FLT_RMT} explicit, and similarly we have for $k=\lfloor n\alpha\rfloor<\min\{n,m\}$,
\begin{align}
	&-n\widetilde g_n(1-\tfrac kn)-n\widetilde \U_n(1)  \label{eq:explicittilde} \\&\qquad=
	k\log\left(\frac{k}{\sigma^2 n}\right)+k
	+(n-k)\log \left(1-\frac{k}{n}\right)	+(m-k)\log \left(1-\frac{k}{m}\right).\nonumber 
\end{align} 
\end{theorem} 
See Figure \ref{fig:sim} for finite-$n$ realizations of the process $\Gau$ for a Wishart matrix $\mathbf Y_n$. Recall that the topology on $\mathcal D((0,\min(1,1/c)))$ fixes a compact subinterval of $(0,\min(1,1/c))$ first, such that the left hand side of \eqref{eq:FLT_RMT} is eventually almost surely finite and well defined and we may set the process equal to zero whenever it is undefined.

The centering law ${\sf MP}_{\sigma^2,n/m}$ may be replaced by $\mathbb E\mu_n$ leading to the same limit. This can be shown via a variance bound from the Efron--Stein inequality, which can be used to derive locally uniform convergence $ \mathbb E \log p_n-n\widetilde \U_n \to0$.

If additionally $m(n)-n/c$ converges, we may also center by the limiting law ${\sf MP}_{\sigma^2,c}$ (corresponding to $\Sti_\infty$  and $g$), yielding a deterministic shift to $\Gau(\alpha)$. If $m(n)-n/c$ does not converge, then such a centering will not work, see for instance the discussion in \cite[\S9.7]{BS10spectral}. We omit the details. Instead, we explicitly present the following implication due to recentering by the expected coefficients.

\begin{corollary}\label{cor:RMT}
With the assumptions and process $\Gau$ of Theorem \ref{thm:FLT_RMT}, we have
	$$
	\left(\frac{e_{\lfloor \alpha n\rfloor:n}}{\mathbb E e_{\lfloor\alpha n\rfloor:n}}\right)_{\alpha\in (0,\min(1,1/c))}\toweak \left(\sqrt{1-c\alpha^2}e^{\Gau(\alpha)} \right)_{\alpha\in (0,\min(1,1/c))}
	$$
	i.e. the normalized coefficient converge to a lognormal process with mean one. 
	\end{corollary}

\begin{remark}[Related literature]\label{rem:literature-RMT}
First, it is known that the expected characteristic polynomial is a rescaled Laguerre polynomial given by \begin{align}\label{eq:ELaguerre}
\E [p_n(-z)]=\Big(\frac{\sigma^2}{m}\Big)^{n}n!L_n^{(m-n)}(mz/\sigma^2),\end{align} 
see for instance \cite{FG06,DumitriuEdelmann,Szego}. 
 It has been shown in \cite{JalKabMarII}, Theorem 6.1 and Corollary 4.14, that also empirical measures of zeros of $\E [p_n(-z)]$ converge weakly to the same ${\sf MP}_{\sigma^2,c}$ and that the profile of the coefficients of $\E [p_n( z)]$ have asymptotics following $\eee^{-g^{\prime}(\alpha)} =\frac{\sigma^2\alpha(c\alpha+1-c)}{1-\alpha}$. This shows that the profile asymptotic holds for expected coefficients and almost surely for the random coefficients. Theorem \ref{thm:FLT_RMT} above reveals the profile fluctuations of those random coefficients. 
	
	Moreover, note that the extreme case of the first finitely many leading coefficients can in principle be deduced from existing results in the literature as follows, see also \cite[Remark 6.2]{JalKabMarII}. The CLT of smooth linear spectral statistics \cite[\S 9.7]{BS10spectral} applies here, which we will also use during the proof. 
	For the test function $\log(1+y)$, the CLT gives  
\begin{align*}
 	&	\log a_{n:n}+n\int\log(1+y)\text{\normalfont\sffamily MP}_{\sigma^2,n/m}(\dd y)\\
 	&\qquad =-\log\det({\bf Y}_n+{\bf Id})+n\int\log(1+y)\text{\normalfont\sffamily MP}_{\sigma^2,n/m}(\mathrm dy)\toweak  - \intprocleftup (1). 
\end{align*}
More generally, for fixed $j\in\N$, the first finitely many coefficients can be expressed as $a_{n-j:n}={a_{n:n}}\cdot {e_{j:n} }$, and via Newton's formula as polynomials in the trace powers $\mathrm{Tr}({\bf Y}_n^1),\dots,\mathrm{Tr}({\bf Y}_n^j)$ all of which again are smooth linear spectral statistics. However, for $j\asymp n$, it is unclear how this approach could lead to such simple exponential profiles and their fluctuations.

The limiting process in Theorem \ref{thm:FLT_RMT} has the same distribution as $(\xi(\sqrt c\alpha))_{\alpha\in (0,\min(1,1/c))}$, where $\xi$ is the following Gaussian analytic function
$$
\xi(z)=\sum_{k\ge1}\frac{Z_k}{\sqrt k} z^k,
$$
with i.i.d. $Z_k\sim\mathcal N(0,1)$. In \cite{LambertPaquette}, 
	it is shown that the Gaussian-$\beta$-characteristic polynomial $\varphi_n$ satisfies $\sqrt{\beta/2}(\log \varphi_n(z)-c_n(z))\toweak \xi(J(z))$, where $c_n$ is some centering and $J$ is the inverse Joukowsky map that we will encounter in the proof as well, see \eqref{eq:q-Joukowski}.	Remarkably, for the coefficient fluctuations in Theorem \ref{thm:FLT_RMT}, the natural parametrization by the coefficient index $\alpha$ leads directly to the simple logarithmic covariance kernel $-\log(1-c\alpha\beta)$, without the need for an additional coordinate transformation. The underlying relation between the spectral and coefficient parametrizations is made explicit in the proof of Theorem \ref{thm:FLT_RMT}, see \eqref{eq:q-Joukowski} and the argument following \eqref{eq:nearlytheend}.
	 
	In \cite{Naj}, the secular coefficients of the circular $\beta$-ensemble are related to the Taylor coefficients of $\exp(\sqrt{2/\beta} \xi_{\C})$, where $\xi_{\C}$ is a complex version of the above.
The distributional boundary value of $\exp(\sqrt{2/\beta}\xi_{\C})$ defines the holomorphic multiplicative chaos. On the other hand, taking the real part of $\xi_\C$, the limit in Corollary \ref{cor:RMT} corresponds to the normalized modulus of this analytic exponential. Similar limiting processes appeared for log-characteristic polynomials of covariance matrices such as \cite[page 112]{johnstone2020testing} and \cite[Proposition 1.2]{dumitriupaquette} and \cite[Corollary 4.5]{HLN08}.  
	 Since $\log p_n$ contains all information on the zeros, this reveals once more the connection that has been established in \cite{JalKabMar:I} between the zeros and coefficients (dealt with in Theorem \ref{thm:FLT_RMT}), even on the level of fluctuations.
\end{remark}

\begin{proof}[Proof of Theorem \ref{thm:FLT_RMT}]
	The claim will follow from Corollary \ref{cor:main_symm}. We shall first verify that \ref{ass:A3} holds with $c_n=n$,	deterministic centering $\widetilde \U_n$ defined in \eqref{eq:MPlogpot}, and a non-degenerate centered Gaussian process $(\intprocleftup (z))_{z\in\mathbb H}$ to be defined below. 
	
	\textbf{Finite dimensional convergence.} We apply \cite[Theorem~9.10 (3)]{BS10spectral} with their $H_n\equiv\delta_{\sigma^2}$ and $y_n=p/n\to y$ corresponding to our $n/m\to c$: For every  function $f$ that is real on $\R$ and analytic on a neighborhood of $[0,c_+]$, defined in \eqref{eq:c+}, we have a CLT for the linear spectral statistic
\begin{align}\label{eq:CLT_LSS}
	n\left(	\int_0^\infty f(x)\mu_n({\rm d}x)	-\int_0^\infty f(x)	{\sf MP}_{\sigma^2,n/m}({\rm d}x)	\right)
	\toweak\mathcal N(0,V(f)),
\end{align} 
	with variance
	\begin{align}\label{eq:VarWish}
		V(f):=\frac{1}{(2\pi)^2}\int_{c_-}^{c_+}\int_{c_-}^{c_+}
		&\left(\frac{f(x)-f(y)}{x-y}\right)^2\nonumber\\
		&\times\frac{4\sigma^4c-(x-c_\ast)(y-c_\ast)}{\sqrt{4\sigma^4c-(x-c_\ast)^2}\sqrt{4\sigma^4c-(y-c_\ast)^2}}	\dd x\dd y,
	\end{align}
	where $c_\ast:=(c_-+c_+)/2=\sigma^2(1+c)$.
	This representation of the variance follows from
	\cite[eq.~(4.28) and Remark~4.2]{LP09CLT}, including the
	factor $1/2$ for complex entries. Choosing $f(x)=\sum_{i=1}^{\ell}\frac{\gamma_i}{z_i+x}$ as in the Cram\'er--Wold device \eqref{eq:A2_fd}, this proves finite-dimensional convergence 	$n(\Sti_n-\widetilde\Sti_n)\to \intprocleft'$ to some Gaussian $\intprocleft'$, where
\begin{multline}
	\widetilde\Sti_n(z):= \int_0^\infty\frac{1}{z+x}{\sf MP}_{\sigma^2,n/m}({\rm d}x)\\
	\ton \Sti_\infty(z):=\int_0^\infty\frac{1}{z+x}{\sf MP}_{\sigma^2,c}({\rm d}x)=\frac{\sqrt{(z+c_-)(z+c_+)}-z+\sigma^2(c-1)}{2c\sigma^2z}.\label{eq:MPStieltjes}
\end{multline} 

	\textbf{Analytic tightness.} We continue to follow the arguments of \cite[\S 9.8]{BS10spectral} and will make use of \cite[Lemma~9.11]{BS10spectral}, which provides contour tightness of the regularized, truncated resolvent process. In order to lift this to tightness of $n(\Sti_n(z)-\widetilde\Sti_n(z))$ in $\mathcal A(\mathbb H)$, we use the contour representation and uniform error estimates therein.\footnote{Alternative proofs based on moment estimates can be deduced from \cite[Propositions~6.4(i), 7.1 and~7.2]{najim2016gaussian}, in particular for the Gaussian-Wishart-model, or \cite[Proposition~4 and (2.20)]{Sh11CLT} for real entries. Also, it is possible to show \ref{ass:A3} directly by the same contour arguments applied to the logarithmic kernel, for which we would need to repeat more steps of \cite{BS10spectral}.}

Fix compact $K\subset\mathbb H$ and set $f_z(x)=1/(z+x)$ for $z\in K$ as above. 
Choose a positively oriented (rectangular) contour $\Gamma$ as in \cite[Section~9.8]{BS10spectral}, with left edge $x_\ell>-\min_{z\in K}\Re z$. With probability tending to one, the support of ${\sf MP}_{\sigma^2,n/m}$ and the spectrum of $\mathbf Y_n$ lie inside $\Gamma$ at a fixed positive distance from it.
Then every $f_z$ is analytic on and inside $\Gamma$, and
$$
 C:=\sup_{z\in K, w\in\Gamma}|f_z(w)|<\infty,\qquad C':=\sup_{z\in K, x\ge0}|f_z'(x)|<\infty.
$$ 
The entry comparison in \cite[Section~9.7.1]{BS10spectral} depends on $C'$ (therein denoted by $K_j$) and is thus  uniform in $z\in K$. On the other hand, $C$ (called $K$ on \cite[p.265]{BS10spectral}) makes the contour error estimate following \cite[Eqn. (9.8.5)]{BS10spectral} uniform in $z\in K$. Note that after truncation of the matrix entries, we can explicitly relate both notations by $M_n(-z) = {n} \bigl(\Sti_n^{\rm truncated}(z)-\widetilde\Sti_n(z)\bigr)$, both centered by a transform of ${\sf MP}_{\sigma^2,n/m}$. The regularized Stieltjes-transform process $\widehat M_n$ is introduced in \cite[Eq.~(9.8.2)]{BS10spectral}, extended by conjugation to the lower half of $\Gamma$, and gives 
\begin{align*}
	n\bigl(\Sti_n(z)-\widetilde\Sti_n(z)\bigr)  =M_n(-z)+o_{\mathbb P}(1) =-\frac{1}{2\pi i}\int_\Gamma f_z(w)\widehat M_n(w)\,\mathrm dw
	+o_{\mathbb P}(1),
\end{align*}
uniformly on $K$. It is now time to call 9.11, which is the lemma on p.~264 of \cite{BS10spectral} establishing tightness of $\widehat M_n$ in the space of continuous functions on the upper half $\Gamma^+$ of the contour $\Gamma$.
 Hence we obtain boundedness in probability, more precisely
$$
\sup_{z\in K}\bigl|n(\Sti_n(z)-\widetilde\Sti_n(z))\bigr|
\le \frac{C \operatorname{length}(\Gamma)}{2\pi}
\|\widehat M_n\|_{\infty,\Gamma^+}+o_{\mathbb P}(1)
=O_{\mathbb P}(1).
$$
Since $K$ was arbitrary, Montel's theorem gives tightness in
$\mathcal A(\mathbb H)$.

\textbf{Verification of \ref{ass:A3}}. So far, we have checked that \ref{ass:A2'} holds, by invoking Proposition \ref{thm:Vitali-Porter-random} similar to \ref{ass:A2''}. Hence, by Proposition \ref{prop:Billing} it remains to verify \eqref{eq:B2}. Apply the scalar CLT \eqref{eq:CLT_LSS} with
	$f(x)=\log(u+x)$ for fixed $u>0$. 
	By the mean value theorem for differentiable functions
$$
	\left|\frac{\log(u+x)-\log(u+y)}{x-y}\right|\le \frac1u,\qquad  x,y\in [c_-,c_+],
$$
hence \eqref{eq:VarWish} gives
\begin{align}  \label{eq:Vleusquared}
	V(\log(u+\cdot))\le \frac{V(\mathrm{id})}{u^2}
	=\frac{\sigma^4c}{u^2}.
	\end{align}
	By the scalar CLT and Chebyshev's inequality for its Gaussian limit,
	\begin{align*}
			\limsup_{n\to\infty}\mathbb P\left(
			\left|\log p_n(u)-n\widetilde \U_n(u)\right|>\varepsilon
			\right)
			&\le \frac{V(\log(u+\cdot))}{\varepsilon^2}
			\le \frac{\sigma^4c}{\varepsilon^2u^2},
			\longrightarrow 0
		\end{align*}
	as $u\to\infty$. Proposition~\ref{prop:Billing} now gives Assumption~\ref{ass:A3}. 
	
\textbf{Covariance Calculation.}
	It remains to derive the claimed covariance of $\intprocleftup$,
	which we transfer to $\Gau(\alpha)=\intprocleftup (\eee^{-g'(1-\alpha)})$
	at the very end of the proof.
	Using polarization in \eqref{eq:VarWish} with
	$f(x)=\log(u+x)+\log(v+x)$ and differentiating under
	the integral, we obtain, for $u,v>0$,
	\begin{align}
	\frac{\dd}{ \dd u}	\frac{\dd}{ \dd v}
		{\rm Cov}[\intprocleftup (u),\intprocleftup (v)]
		&=\frac{1}{(2\pi)^2}
		\int_{c_-}^{c_+}\int_{c_-}^{c_+}
		\frac{1}{(u+x)(v+x)(u+y)(v+y)}
		\nonumber\\
		&\qquad\times
		\frac{4\sigma^4c-(x-c_\ast)(y-c_\ast)}
		{\sqrt{4\sigma^4c-(x-c_\ast)^2}
			\sqrt{4\sigma^4c-(y-c_\ast)^2}}
		\dd x\dd y \label{eq:Wishart_limiting_field} \\
		&=	\frac14 h_+(u,v)h_-(u,v)\label{eq:hproduct}.
	\end{align}  
In the last line, the double integral has been replaced by a product of
	\begin{align}\label{eq:h}
		h_\pm(u,v)
		:=\frac1\pi\int_{c_-}^{c_+}
		\frac{2\sigma^2\sqrt c\pm(x-c_\ast)}
		{(u+x)(v+x)}
		\frac{{\rm d}x}
		{\sqrt{4\sigma^4c-(x-c_\ast)^2}},
	\end{align}
whose cross terms indeed cancel out by symmetry. In order to evaluate $h_\pm$, we first rearrange the integrand into
	\begin{align}\label{eq:split}
		\frac{2\sigma^2\sqrt c\pm(x-c_\ast)}
		{(u+x)(v+x)}
		&=
		\frac{2\sigma^2\sqrt c\mp(u+c_\ast)}{v-u}
		\left(\frac1{u+x}-\frac1{v+x}\right)
		\pm\frac1{v+x}.
	\end{align}
	Consequently, we are left with integrals like
	\begin{align}\label{eq:arcsine}
		\frac1\pi\int_{c_-}^{c_+}
		\frac1{u+x}
		\frac{{\rm d}x}
		{\sqrt{4\sigma^4c-(x-c_\ast)^2}}
		=
		\frac1{\sqrt{(u+c_-)(u+c_+)}}.
	\end{align}
	This can either be seen as the Stieltjes transform of the arcsine distribution centered at $c_\ast$, or follows by a standard Beta-integral (for instance \cite[5.12.4]{NIST} with $a=b=1/2$ after a linear transformation). Combining \eqref{eq:h}, \eqref{eq:split} and \eqref{eq:arcsine}, we obtain
	\begin{align} 
		h_\pm(u,v)
	=&\frac {2\sigma^2\sqrt c\mp (u+c_*)}{v-u}\left(	\frac 1 {\sqrt{(u+c_-)(u+c_+)}}-\frac 1 {\sqrt{(v+c_-)(v+c_+)}} \right)\pm \frac 1 {\sqrt{(v+c_-)(v+c_+)}}\nonumber \\
	=&\frac 1 {v-u}\left(	\frac{2\sigma^2\sqrt c\mp (u+c_*)}{\sqrt{(u+c_-)(u+c_+)}}-\frac{2\sigma^2\sqrt c\mp (v+c_*)}{\sqrt{(v+c_-)(v+c_+)}} \right). \label{eq:h2}
	\end{align}
	The identities at $u=v$ are understood by continuity. In order to simplify this, recall the explicit representation \eqref{eq:MPStieltjes}  in the more convenient notation
	\begin{align}\label{eq:q-Joukowski}
		q(u):=\sqrt c (1-u\Sti_\infty(u))= \frac{u+c_\ast-\sqrt{(u+c_-)(u+c_+)}}{2\sqrt c\sigma^2},
	\end{align}
	where we used $c_\ast=(c_{-}+c_{+})/2=\sigma ^2(1+c)$. Moreover, we obtain a variant of the well known self consistent (Dyson) equation\footnote{It may also be seen as coordinates of a Joukowsky transform. Note furthermore that these are equivalent \eqref{eq:g'} for $\alpha=\Phi_\infty(u)$ and using $\exp(-g'(\Phi_\infty (u)))=u$.} of the Stieltjes transform of ${\sf MP}_{\sigma^2,c}$ taking the form
	\begin{align}\label{eq:selfconseqn}
		u+c_\ast 
		= \sigma ^2 \sqrt c \left(q(u)+\frac{1}{q(u)}\right),
	\end{align} 
	and the analogous equation
	$$\sqrt{(u+c_-)(u+c_+)}=\sigma ^2\sqrt c\left(\frac 1 {q(u)}-q(u)\right).$$
	Therefore, we found that \eqref{eq:h2} simplifies to
	\begin{align*}
		h_\pm(u,v)=\frac 1 {v-u}\left(\frac{q(u)\mp 1}{1\pm q(u)}-\frac{q(v)\mp 1}{1\pm q(v)} \right)=\frac{ q(u)-q(v)}{v-u}\cdot \frac 2 {(1\pm q(u))(1\pm q(v))} .
	\end{align*}
	Furthermore, we may subtract \eqref{eq:selfconseqn} with $v$ and $u$, giving
	$$v-u= 
	\sigma ^2\sqrt c\big(q(u)-q(v)\big)\frac{1-q(u)q(v)} {q(u)q(v)}.$$
	Hence, we obtain
	\begin{align*}
		h_\pm(u,v)= \frac {2q(u)q(v)}{\sigma^2\sqrt c(1-q(u)q(v))(1\pm q(u))(1\pm q(v))},
	\end{align*}
	which upon multiplication gives the right-hand side of~\eqref{eq:hproduct}:
	\begin{align}
		\frac 1 4 h_+(u,v)h_-(u,v)= \frac {q(u)^2q(v)^2}{\sigma^4 c(1-q(u)q(v))^2(1 -q(u)^2)(1-q(v)^2)} 
	\end{align}
	On the other hand, differentiation of \eqref{eq:selfconseqn} yields
	$$q'(u)=\frac{q(u)^2}{\sigma ^2\sqrt c (q(u)^2-1)}.$$
	Finally, we arrive at the simple identity
	\begin{align}\label{eq:h+h-}
		\frac 1 4 h_+(u,v)h_-(u,v)= \frac { q'(u)q'(v)}{ (1-q(u)q(v))^2 }=- \frac{\dd}{\dd u}\frac{\dd}{\dd v}\log\big( 1-q(u)q(v)\big).
	\end{align} 
By \eqref{eq:Wishart_limiting_field}, it remains to determine the integration constants to be zero, for which we recall that \eqref{eq:Vleusquared} gives ${\rm Var}[\intprocleftup (u)] =V(\log(u+\cdot))\lesssim u^{-2}$. By Cauchy--Schwarz, ${\rm Cov}[\intprocleftup (u),\intprocleftup (v)]\to0$ for $u\to\infty$ or $v\to\infty$.
	Together with $q(\infty)=0$, integration of
	\eqref{eq:hproduct} and \eqref{eq:h+h-} therefore gives
	\[
	{\rm Cov}[\intprocleftup (u),\intprocleftup (v)]
	=-\log(1-q(u)q(v)),
	\qquad u,v>0.
	\]
	By analytic continuation,  we obtain both the covariance and pseudo-covariance kernels
		\begin{align}\label{eq:nearlytheend}
		\begin{cases}
		\mathbb{E}[\intprocleftup (z_1)\overline{\intprocleftup (z_2)}]
		=
		-\log\bigl(1-q(z_1)\overline{q(z_2)}\bigr),\\
		\mathbb{E}[\intprocleftup (z_1)\intprocleftup (z_2)]
		=
		-\log\bigl(1-q(z_1)q(z_2)\bigr),
		\end{cases}
		\qquad z_1,z_2\in\mathbb H.
	\end{align}
	Together, these kernels determine the law of the centered complex Gaussian process.
	As anticipated, the statement now follows from Corollary \ref{cor:main_symm} and setting 	$\Gau(\alpha)=\intprocleftup (\eee^{-g'(1-\alpha)})$. Recall that \eqref{eq:phi_Phi_inversion} implies $x\Sti_\infty(x)=1-\alpha$ for $x=\eee^{-g'(1-\alpha)}$, so $q(x)=\sqrt c\alpha$ by \eqref{eq:q-Joukowski} yields the claimed covariance
	\begin{align*}
		{\rm Cov}[\Gau(\alpha),\Gau(\beta)]
		&=-\log(1-c\alpha\beta).
	\end{align*}
\textbf{Explicit centering.} It remains to verify \eqref{eq:explicittilde}.
Since $\widetilde \U_n$ is the logarithmic potential of $\mathsf{MP}_{\sigma^2,n/m}$, the profile $\widetilde g_n$ is the Marchenko--Pastur profile at aspect ratio $c$ replaced by $n/m$. It follows from \eqref{eq:g'} that
\begin{align}\label{eq:irgendeineformel}
	\frac{\mathrm d}{\mathrm dt}
	\bigl(\widetilde g_n(1-t)+\widetilde \U_n(1)\bigr)
	&=-\widetilde g_n'(1-t)
	=\log\frac{\sigma^2(1-t)(1-(n/m)t)}{t}.
\end{align}
Note that by \eqref{eq:tildeg_n_def_as_Legendre2} and the normalization of $\widetilde U_n$ at infinity,
$$
\widetilde g_n(1)+\widetilde U_n(1)
=\lim_{x\to\infty}\bigl[\widetilde U_n(x)-\log x+(1-x\widetilde U_n'(x))\log x\bigr]=0,
$$
since $0\le 1-x\widetilde U_n'(x)\lesssim 1/x$. Integrating \eqref{eq:irgendeineformel} from $t=0$ to $k/n$, then yields
\begin{align*}
	n\bigl(\widetilde g_n(1-k/n)+\widetilde \U_n(1)\bigr)
	&=n\int_0^{k/n}\log\frac{\sigma^2(1-t)(1-(n/m)t)}{t}\,\mathrm dt\\
	&=k\log\frac{\sigma^2 n}{k}-k
	-(n-k)\log\!\left(1-\frac{k}{n}\right)
	-(m-k)\log\!\left(1-\frac{k}{m}\right)
\end{align*}
as claimed in \eqref{eq:explicittilde} after changing the sign.  Alternatively,~\eqref{eq:explicittilde} can be derived using the explicit formula for the logarithmic potential of the Marchenko-Pastur distribution given in \cite[Eq.~(37)]{JalKabMarII}.
\end{proof}

\begin{proof}[Proof of Corollary \ref{cor:RMT}] Extracting the coefficient from the Laguerre identity \eqref{eq:ELaguerre} gives
$$
\E e_{k:n}=\binom nk\frac{m!}{(m-k)!}
\left(\frac{\sigma^2}{m}\right)^k .
$$
For $k=\lfloor n\alpha\rfloor$, Stirling’s formula gives locally uniformly in $\alpha\in (0,\min(1,1/c)$,
\begin{align*}
	\log\mathbb E e_{k:n}
	={}&k\log\frac{\sigma^2 n}{k}-k
	-(n-k)\log\!\left(1-\frac{k}{n}\right)
	-(m-k)\log\!\left(1-\frac{k}{m}\right)\\
	&-\frac12\log\!\left(2\pi k\left(1-\frac{k}{n}\right)
	\left(1-\frac{k}{m}\right)\right)+o(1).
\end{align*} 
Compared to the deterministic centering in \eqref{eq:FLT_RMT} represented as in \eqref{eq:explicittilde} and \eqref{eq:g''}, the above exceeds it by $-\log\sqrt{1-c\alpha^2}+o(1)$, which yields the claimed limit after exponentiating.
\end{proof}

\section{Proof of Theorem~\ref{thm:main}}\label{sec:proof}
In this section, we will make use of the probabilistic representations \eqref{eq:coeff} and \eqref{eq:P_n}. We will stick to the notation $\Phi_n$, which is most convenient for inversions.
\subsection{Preliminaries}
Note that if $f\in\mathcal{A}(\mathbb{H})$, then
$z\mapsto \int_{1}^z f(u)  {\rm d}u$ is analytic on $\mathbb{H}$ and 
the mapping\footnote{ The contour of integration is taken inside $\mathbb{H}$.}
$$
\mathcal{A}(\mathbb{H})\ni f~\longmapsto~\int_1^{(\cdot)} {f(u)} {\rm d}u\in\mathcal{A}(\mathbb{H})
$$
is continuous with respect to the locally uniform topology on $\mathcal{A}(\mathbb{H})$, and the same holds for differentiation.
Recalling \eqref{eq:Phi_Psi_integral}, the continuous mapping theorem immediately yields

\begin{lemma}\label{lem:phi_psi_conv}
If \ref{ass:A1} and \ref{ass:A2'} hold, then
\begin{equation}\label{eq:psi_phi_joint_n_convergence_hp}
\left(c_n\left(\Psi_n(z)-\widetilde\Psi_n(z)\right),c_n(\Phi_n(z)-\widetilde\Phi_n(z))\right)_{z\in\mathbb{H}}~\toweak~\big(\intprocleft(z),z\intprocleft'(z)\big)_{z\in\mathbb{H}},
\end{equation}
in the space of probability measures on $\mathcal{A}(\mathbb{H})\times\mathcal{A}(\mathbb{H})$. In particular,
\begin{equation}\label{eq:psi_phi_joint_n_convergence}
\left(c_n\left(\Psi_n(x)-\widetilde\Psi_n(x)\right),c_n(\Phi_n(x)-\widetilde\Phi_n(x))\right)_{x>0}~\toweak~\left(\intprocleft(x),x\intprocleft'(x)\right)_{x>0},
\end{equation}
in the space of probability measures on $\mathcal C((0,\infty),\mathbb{R}^2)$. 
\end{lemma}

\subsection{Exponential tilting}

Recall that $\mathcal{F}_n$ denotes the sigma-algebra generated by $\lambda_{1:n},\ldots,\lambda_{n:n}$, for $n\in\mathbb{N}$, and $S_n:=\sum_{j=1}^{n}B_{j:n}$ as in \eqref{eq:P_n}, where $P_n(z)=\mathbb{E}[z^{S_n}|\mathcal{F}_n]$. For any $\theta>0$, introduce a new (tilted) conditional probability measure $\widehat{\mathbb{P}}_{\theta}[\cdot | \mathcal{F}_n]$ (with the corresponding conditional expectation $\widehat{\mathbb{E}}_{\theta}[\cdot  | \mathcal{F}_n]$) by the identity
\begin{equation}\label{eq:exponential_tilting}
\widehat{\mathbb{E}}_{\theta}[f(B_{1:n},\ldots,B_{n:n}) | \mathcal{F}_n]=\frac{\mathbb{E}[f(B_{1:n},\ldots,B_{n:n})\theta^{S_n} | \mathcal{F}_n]}{P_n(\theta)},
\end{equation}
for any bounded measurable $f:\mathbb{R}^n\to \mathbb{R}$. Substitution $f(x_1,\ldots,x_n)=\1[x_1+\cdots+x_n=k]$ into~\eqref{eq:exponential_tilting} results in
\begin{equation}\label{eq:exponential_tilting_probabilities1}
a_{k:n}=\mathbb{P}[S_n=k | \mathcal{F}_n]=P_n(\theta)\theta^{-k}\widehat{\mathbb{P}}_{\theta}[S_n=k | \mathcal{F}_n].
\end{equation}
Recalling the definition~\eqref{eq:g_n_def} of $g_n$ and the representation~\eqref{eq:tildeg_n_def_as_Legendre} of $\widetilde g_n$, we obtain 
\begin{align} 
c_n\big (g_n(\alpha)-\widetilde g_n(\lfloor n\alpha \rfloor/n)+\tfrac 1 n  \log \sqrt n\big )
=&-c_n\left(\frac{\lfloor n\alpha \rfloor}{n}\log \theta - \frac{\lfloor n\alpha \rfloor}{n}\log \widetilde\Phi_n^{\leftarrow}(\lfloor n\alpha \rfloor/n)\right)\notag\\
& +c_n\left(\frac{\log P_n(\theta)}{n}-\widetilde\Psi_n(\widetilde\Phi_n^{\leftarrow}(\lfloor n\alpha \rfloor/n))\right)\notag\\
&+\frac{c_n}{n}\log\left(\sqrt{n}\cdot \widehat{\mathbb{P}}_{\theta}[S_n=\lfloor n\alpha \rfloor | \mathcal{F}_n]\right) \label{eq:p_vs_p_theta} 
\end{align}
 for $\alpha \in [0,1]$. The idea is to pick an $\mathcal{F}_n$-measurable $\theta_{\ast}=\theta_{\ast}(\alpha ,n)$ such that $\widehat{\mathbb{E}}_{\theta_{\ast}}[S_n | \mathcal{F}_n]= \lfloor n\alpha \rfloor$. As we shall check below, see Lemma~\ref{lem:third_summand}, this choice guarantees that $\sqrt{n}\cdot \widehat{\mathbb{P}}_{\theta_{\ast}}[S_n=\lfloor n\alpha \rfloor | \mathcal{F}_n]$ converges a.s. locally uniformly on $(\underline{m},\overline{m})$ to a finite deterministic constant. Therefore, the last summand in~\eqref{eq:p_vs_p_theta} may be negligible if $c_n=o(n)$ as claimed in Remark~\ref{rem:scaling}, or gives the term as claimed in \eqref{eq:main_flt_measures}. The principal contribution comes from the first two summands in~\eqref{eq:p_vs_p_theta}.

By choosing $f(x_1,\ldots,x_n)=x_1+\cdots+x_n$ in~\eqref{eq:exponential_tilting} we obtain the equation
\begin{equation*}
\widehat{\mathbb{E}}_{\theta_{\ast}}[S_n | \mathcal{F}_n]=\frac{\mathbb{E}[S_n\theta_{\ast}^{S_n} | \mathcal{F}_n]}{P_n(\theta_{\ast})}=\lfloor n\alpha \rfloor~\Longleftrightarrow~\frac{\theta_{\ast} P^{\prime}_n(\theta_{\ast})}{n P_n(\theta_{\ast})}=\frac{\lfloor n\alpha \rfloor}{n},\quad \alpha \in (0,1),
\end{equation*}
which takes the form $\Phi_n(\theta_{\ast})=\frac{\lfloor n\alpha \rfloor}{n}$, where $\Phi_n$ was defined in~\eqref{eq:phi_n_def_and_conv}.

The equation $\Phi_n(\theta_{\ast})=\frac{\lfloor \alpha n \rfloor}{n}$ has an a.s.~unique solution 
\begin{align}
\theta_{\ast}:=\theta_{\ast}(\alpha ,n):=\Phi_n^{\leftarrow}(\tfrac{\lfloor \alpha n \rfloor}{n}),\quad \alpha \in (\underline{m}, \overline{m}), \label{eq:theta*}
\end{align} 
for $n$ sufficiently large, since the function $x\mapsto \Phi_n(x)$ is eventually strictly increasing and continuous on $(0,\infty)$. Assumption \ref{ass:A1'} yields
\begin{equation}\label{eq:JalKabMar2}
\lim_{n\to\infty}\sup_{\alpha \in K}|\Phi_n^{\leftarrow}(\alpha )-\Phi_{\infty}^{\leftarrow}(\alpha )|=0\qquad\text{a.s.~for every compact }K\subset (\underline{m},\overline{m}),
\end{equation}
hence, 
\begin{equation}
\lim_{n\to\infty}\sup_{\alpha \in K}|\theta_{\ast}(\alpha ,n)-\Phi^{\leftarrow}_{\infty}(\alpha )|=0\qquad\text{a.s.~for every compact }K\subset (\underline{m},\overline{m}).\label{eq:JalKabMar3}
\end{equation}

\subsection{Controlling three terms}
In order to prove Theorem~\ref{thm:main}, we have to deal with all three terms of \eqref{eq:p_vs_p_theta} and we begin with the convergence of the last term.
\begin{lemma}\label{lem:third_summand}
	Fix a compact set $K\subset (\underline{m},\overline{m})$ and let $\theta_{\ast}(\alpha ,n)=\Phi_n^{\leftarrow}(\lfloor n\alpha \rfloor/n)$. Then, as $n\to\infty$,
	$$
	\sup_{\alpha \in K}\left|\sqrt{n}\cdot \widehat{\mathbb{P}}_{\theta_{\ast}(\alpha ,n)}[S_n=\lfloor n\alpha \rfloor | \mathcal{F}_n]-\sqrt{\frac{-g''(\alpha )}{2\pi}}\right|\to 0\quad\text{a.s.}
	$$
\end{lemma}
\begin{proof}
Assume that there exist an event of positive probability, on which there are sequences $(n_m)_{m\in\mathbb{N}}\subseteq\mathbb{N}$ and $(\alpha _m)_{m\in\mathbb{N}}\subseteq K$ such that
	\begin{align}\label{eq:tilted_probability_sqrt_n_asymp}
		\lim_{m\to\infty}\Big|\sqrt{n_m}\cdot \widehat{\mathbb{P}}_{\theta_{\ast}(\alpha _m,n_m)}[S_{n_m}=\lfloor n_m \alpha _m\rfloor | \mathcal{F}_{n_m}]
		-\sqrt{\frac{-g''(\alpha_m)}{2\pi}}\Big|>0 .
	\end{align}
	Passing to further subsequence we can and do assume that $\alpha _m\to \alpha \in (\underline{m},\overline{m})$.
	
	On the other hand, according to Eq.~(50) in~\cite{JalKabMar:I}, condition \ref{ass:A1} implies that a.s.
	\begin{equation}\label{eq:tilted_probability_sqrt_n_asymp1}
		\widehat{\mathbb{P}}_{\theta_{\ast}(\alpha _m,n_m)}[S_{n_m}=\lfloor n_m \alpha _m\rfloor | \mathcal{F}_{n_m}]~\sim~\frac{1}{\sqrt{2\pi \sigma_{n_m}^2(\alpha _m)}},\quad m\to\infty,
	\end{equation}
	where
	$$
	\sigma_n^2(\alpha ):=\widehat{{\rm Var}}_{\theta_{\ast}(\alpha ,n)}[S_n|\mathcal{F}_n]=\sum_{j=1}^{n}\widehat{\mathbb{P}}_{\theta_{\ast}(\alpha ,n)}[B_{j:n}=1|\mathcal{F}_n]\widehat{\mathbb{P}}_{\theta_{\ast}(\alpha ,n)}[B_{j:n}=0|\mathcal{F}_n].
	$$
	Using~\eqref{eq:exponential_tilting} with $f(x_1,\ldots,x_n)=\1[x_j=1]$ and recalling $\eta_{j:n}=1/(1+\lambda_{j:n})$ from \eqref{eq:eta_j}, we conclude that
	$$
	\sigma_n^2(\alpha )=\sum_{j=1}^{n}\frac{\theta_{\ast}(\alpha ,n)\eta_{j:n}(1-\eta_{j:n})}{(1-\eta_{j:n}+\eta_{j:n}\theta_{\ast}(\alpha ,n))^2}=n\theta_{\ast}(\alpha ,n)\int_0^\infty\frac{x}{(x+\theta_{\ast}(\alpha ,n))^2}\mu_n({\rm d}x).
	$$
	By P\'olya's extension of Dini's theorem
	$$
	\lim_{n\to\infty}\sup_{z\in K_1}\left|\int_0^\infty\frac{x}{(x+z)^2}\mu_n({\rm d}x)-\int_0^\infty\frac{x}{(x+z)^2}\mu_{\infty}({\rm d}x)\right|=0\qquad\text{a.s.},
	$$
	where $K_1\subset (0,\infty)$ is any fixed compact. Therefore,~\eqref{eq:JalKabMar3} implies
	$$
	\lim_{n\to\infty}\sup_{\alpha \in K}\left|\frac{\sigma_n^2(\alpha )}{n}-\int_0^\infty\frac{x\Phi_{\infty}^{\leftarrow}(\alpha )}{(x+\Phi_{\infty}^{\leftarrow}(\alpha ))^2}\mu_{\infty}({\rm d}x)\right|=0\qquad\text{a.s.}
	$$
	Plugging this into~\eqref{eq:tilted_probability_sqrt_n_asymp1} and recalling that $\alpha _m\to \alpha $, we obtain almost surely as $m\to\infty$
	\begin{align}\label{eq:onemorestep}
		\sqrt{{n_m}}\cdot\widehat{\mathbb{P}}_{\theta_{\ast}(\alpha _m,n_m)}[S_{n_m}=\lfloor n_m \alpha _m\rfloor | \mathcal{F}_{n_m}] 
		\to~\left(2\pi \int_0^\infty\frac{x\Phi_{\infty}^{\leftarrow}(\alpha )}{(x+\Phi_{\infty}^{\leftarrow}(\alpha ))^2}\mu_{\infty}({\rm d}x)\right)^{-1/2}.
	\end{align}
	A simple calculation shows that 
	\begin{align*}
	\Phi_{\infty}'(z)=\int_0^\infty\frac{x}{(x+z)^2}\mu_{\infty}({\rm d}x),
	\end{align*}
	hence for $z=\Phi_{\infty}^{\leftarrow}(\alpha )$ and differentiating the inverse,
	\begin{align}\label{eq:gtophi}
		\int_0^\infty\frac{x\Phi_{\infty}^{\leftarrow}(\alpha )}{(x+\Phi_{\infty}^{\leftarrow}(\alpha ))^2}\mu_{\infty}({\rm d}x)=\Phi_{\infty}^{\leftarrow}(\alpha )\Phi_{\infty}'(\Phi_{\infty}^{\leftarrow}(\alpha ))=\frac{\Phi_{\infty}^{\leftarrow}(\alpha )}{(\Phi_{\infty}^{\leftarrow})'(\alpha )}.
	\end{align}
	Finally by \eqref{eq:phi_Phi_inversion}, we have $\eee^{-g'(\alpha)}= \Phi_{\infty}^{\leftarrow}(\alpha)$, thus $$(\Phi_{\infty}^{\leftarrow})'(\alpha )=-g''(\alpha )\eee^{-g'(\alpha )}=-g''(\alpha ) \Phi_{\infty}^{\leftarrow}(\alpha ).$$
	Therefore, \eqref{eq:onemorestep} contradicts~\eqref{eq:tilted_probability_sqrt_n_asymp} thereby finishing the proof.
\end{proof}

Now let us turn to the first summand in \eqref{eq:p_vs_p_theta}.
\begin{lemma}\label{lem:first_summand} 
	For $\alpha \in (\underline{m}, \overline{m})$, define
	\begin{align}\label{eq:Y}
		 Y(\alpha ):=	-\frac{ \intprocleft'(\Phi_{\infty}^{\leftarrow}(\alpha ))}{\Phi_{\infty}^{\prime}(\Phi_{\infty}^{\leftarrow}(\alpha ))}.
	\end{align}
Then, the following convergences holds on $\mathcal D((\underline{m}, \overline{m}))$ 
\begin{align}
	c_n\left(\Phi_n^{\leftarrow}(\alpha)-\widetilde\Phi_n^{\leftarrow}(\alpha)\right)&\toweak~ \Phi_{\infty}^{\leftarrow}(\alpha )Y(\alpha),\label{eq:third_summand_prep}
\\
	c_n\left(\frac{\lfloor n\alpha \rfloor}{n}\log \theta_{\ast}(\alpha ,n) - \frac{\lfloor n\alpha \rfloor}{n}\log \widetilde\Phi_n^{\leftarrow}(\lfloor n\alpha \rfloor/n)\right)&\toweak  {\alpha  }  Y(\alpha )\label{eq:third_summand}
	\end{align}

\end{lemma}
Here and in the following, we shall drop the subscripts $(\cdot)_{\alpha \in (\underline{m}, \overline{m})}$ for notational convenience.
\begin{proof} 
 Let $K$ be a fixed compact subset of $(\underline{m},\overline{m})$. Using that both $\widetilde\Phi_n^{\leftarrow}\to \Phi_{\infty}^{\leftarrow}$ and  $\Phi_n^{\leftarrow}\to \Phi_{\infty}^{\leftarrow}$ uniformly on $K$ and continuity of the composition, we conclude from~\ref{ass:A2'} that
\begin{align}
	c_n(\Phi_n(\widetilde\Phi_n^{\leftarrow}(\alpha))-\alpha)_{\alpha \in K}~&\toweak~(\Phi_{\infty}^{\leftarrow}(\alpha )\intprocleft'(\Phi^{\leftarrow}_{\infty}(\alpha)))_{\alpha \in K},\\
	c_n(\widetilde\Phi_n(\Phi^{\leftarrow}_n(\alpha))-\alpha)_{\alpha \in K}
	~&\toweak~(-\Phi_{\infty}^{\leftarrow}(\alpha )\intprocleft'(\Phi^{\leftarrow}_{\infty}(\alpha)))_{\alpha \in K}\label{eq:phi_inverse_converge1}
\end{align}
on $\mathcal C((\underline{m}, \overline{m}))$. In order to prove \eqref{eq:third_summand_prep}, we argue along the lines of Theorem 13.7.2 in~\cite{Whitt:2002}.
By the mean value theorem, there exist some $\beta(n,\alpha)\to\Phi_\infty^{\leftarrow}(\alpha)$ such that
$$
\widetilde\Phi_n(\Phi^{\leftarrow}_n(\alpha))-\alpha=\widetilde\Phi_n(\Phi^{\leftarrow}_n(\alpha))-\widetilde\Phi_n(\widetilde\Phi_n^{\leftarrow}(\alpha))=
\widetilde\Phi_n^{\prime}(\beta(n,\alpha))(\Phi^{\leftarrow}_n(\alpha)-\widetilde\Phi_n^{\leftarrow}(\alpha)).
$$
From~\eqref{eq:centering_to_phi_infty} and analyticity we conclude that $\widetilde\Phi_n^{\prime}~\to~\Phi^{\prime}_{\infty}$ locally uniformly on $\mathbb{H}$.  
Combining this with \eqref{eq:phi_inverse_converge1}, we arrive at~\eqref{eq:third_summand_prep}. 

Recalling the definition \eqref{eq:theta*} of $\theta_{\ast}$, it follows
\begin{equation}\label{eq:phi_inverse_converge2}
c_n(\theta_{\ast}(\alpha ,n)-\widetilde\Phi_n^{\leftarrow}(\lfloor n\alpha \rfloor/n))
=c_n(\Phi_n^{\leftarrow}(\lfloor n\alpha \rfloor/n)-\widetilde\Phi_n^{\leftarrow}(\lfloor n\alpha \rfloor/n))
~\toweak~\Phi_{\infty}^{\leftarrow}(\alpha )Y(\alpha )
\end{equation}
on $\mathcal D((\underline{m}, \overline{m}))$. Note that $\Phi_{\infty}^{\leftarrow}(\alpha)$ does not vanish on $(\underline{m}, \overline{m})$ and hence
$$
\liminf_{n\to\infty}\widetilde\Phi_n^{\leftarrow}(\lfloor n\alpha\rfloor/n)>0.
$$
Therefore,~\eqref{eq:phi_inverse_converge2} implies for every compact $K\subset (\underline{m},\overline{m})$,
$$
c_n \sup_{\alpha  \in K}\left(\frac{\theta_{\ast}(\alpha ,n)}{\widetilde\Phi_n^{\leftarrow}(\lfloor n\alpha \rfloor/n)}-1\right)^2~\ton~0 $$
in probability.
Hence, by Taylor's expansion of $x\mapsto \log (1+x)$,
\begin{align}\label{eq:Log_of_theta}
c_n
\left(\log \theta_{\ast}(\alpha ,n)-\log \widetilde\Phi_n^{\leftarrow}(\lfloor n\alpha \rfloor/n)\right)
= c_n\log \left(1+\left(\frac{\theta_{\ast}(\alpha ,n)}{\widetilde\Phi_n^{\leftarrow}(\lfloor n\alpha \rfloor/n)}-1\right)\right)
~\toweak~  {Y(\alpha )} 
\end{align}
in the space of probability measures on $\mathcal D((\underline{m}, \overline{m}))$.  Summarizing, we have shown that the first summand in~\eqref{eq:p_vs_p_theta} with $\theta=\theta_{\ast}(\alpha ,n)$ converges weakly on the space $\mathcal D((\underline{m}, \overline{m}))$ to the process $(\alpha Y(\alpha ) )_{\alpha \in (\underline{m}, \overline{m})}$ as claimed.
\end{proof}

For subsequent usage we also note that the same argument implies that
\begin{equation}\label{eq:Psi_of_theta}
	c_n\big (\widetilde\Psi_n(\theta_{\ast}(\alpha ,n))-\widetilde\Psi_n(\widetilde\Phi_n^{\leftarrow}(\lfloor n\alpha \rfloor/n))\big )
	~\toweak \Psi_{\infty}^{\prime}(\Phi_{\infty}^{\leftarrow}(\alpha ))\Phi_{\infty}^{\leftarrow}(\alpha )Y(\alpha )
\end{equation}
on $\mathcal D((\underline{m}, \overline{m}))$.
For the derivation of both~\eqref{eq:Log_of_theta} and~\eqref{eq:Psi_of_theta} we actually also used a well known variation of the `delta method'. Let us deal with the second summand in~\eqref{eq:p_vs_p_theta}.
\begin{lemma}\label{lem:second_summand}
	Jointly with \eqref{eq:third_summand}, we have on  $\mathcal D((\underline{m}, \overline{m}))$,
	\begin{align}\label{eq:log_potential_converge3}
		c_n\left(\frac{1}{n}\log P_n(\theta_{\ast}(\alpha ,n))-\widetilde\Psi_n (\widetilde\Phi_n^{\leftarrow}(\tfrac{\lfloor \alpha n \rfloor}n))\right)
		~\toweak~\intprocleft(\Phi_{\infty}^{\leftarrow}(\alpha ))+\Psi_{\infty}^{\prime}(\Phi_{\infty}^{\leftarrow}(\alpha ))\Phi_{\infty}^{\leftarrow}(\alpha )Y(\alpha ) .
	\end{align}
\end{lemma}

\begin{proof} By \eqref{eq:psi_def}, we write $\frac 1 n \log P_n=\Psi_n$ and by \eqref{eq:psi_phi_joint_n_convergence}, we have $c_n(\Psi_n(z)-\widetilde\Psi_n(z)) \toweak~ \intprocleft(z)$, which gives the first term for $z=\theta_{\ast}(\alpha ,n)\to \Phi_{\infty}^{\leftarrow}(\alpha )$ due to \eqref{eq:JalKabMar3}. The remaining term is precisely~\eqref{eq:Psi_of_theta} giving the second term in the claim. The claim follows from the fact that the two distributional convergences hold jointly by Lemma~\ref{lem:phi_psi_conv}.
\end{proof}

\subsection{Proof of Theorem~\ref{thm:main}}
By Lemma \ref{lem:first_summand} and Lemma \ref{lem:second_summand}, the sum of the first two terms on the right-hand side of~\eqref{eq:p_vs_p_theta} converges weakly to
$$
- {\alpha Y(\alpha )} + \intprocleft(\Phi_{\infty}^{\leftarrow}(\alpha ))+\Psi_{\infty}^{\prime}(\Phi_{\infty}^{\leftarrow}(\alpha ))\Phi_{\infty}^{\leftarrow}(\alpha )Y(\alpha )
= \intprocleft(\Phi_{\infty}^{\leftarrow}(\alpha ))
$$
on $\mathcal D((\underline{m},\overline{m}))$. Here, we have used $\Phi_{\infty}(z)=z\Psi_{\infty}^{\prime}(z)$ to derive the last equality. By Lemma~\ref{lem:third_summand}, we obtain another term $\log\sqrt{\frac{-g''(\alpha )}{2\pi}}$ uniformly in $\alpha $, which is on the left hand side of our claim \eqref{eq:main_flt_measures}. The proof of Theorem~\ref{thm:main} is finished. 
\qedhere

We still owe the reader the
\begin{proof}[Proof of Proposition \ref{prop:Billing}]
The essential idea is that condition~\eqref{eq:B2} recovers the random additive constant lost when normalizing the logarithmic potential at $1$.

\noindent
Assume \ref{ass:A2} and \eqref{eq:B2}. Fix $x>1$ and recall that we also assume $\widetilde\Psi_n(z)=\widetilde \U_n(z)-\widetilde \U_n(1)$. The continuous mapping theorem yields
\begin{align*} 
&	\left(c_n\left(\U_n(z)-\widetilde \U_n(z)\right)
		-
		c_n\left(\U_n(x)-\widetilde \U_n(x)\right)\right)_{z\in\mathbb{H}}\\
		&\qquad=\left(c_n\bigl[(\Psi_n(z)-\widetilde\Psi_n(z))
		-(\Psi_n(x)-\widetilde\Psi_n(x))\bigr]\right)_{z\in\mathbb{H}}\\
%
	&\qquad\toweak
	\left(\intprocleft(z)-\intprocleft(x)\right)_{z\in\mathbb H}.
\end{align*} 
We now show that $\intprocleft(x)$ converges in probability as $x\to+\infty$. For any fixed $x,y>1$, assumption \ref{ass:A2} and the Portmanteau theorem imply, for every $\varepsilon>0$,
\begin{align*}
\mathbb{P}[|\intprocleft(x)-\intprocleft(y)|>\varepsilon]&\leq
\liminf_{n\to\infty}
\mathbb{P}[|c_n\left(\U_n(x)-\widetilde \U_n(x)\right)
		-
		c_n\left(\U_n(y)-\widetilde \U_n(y)\right)|>\varepsilon]\\
&\leq \limsup_{n\to\infty}\mathbb{P}[|c_n\left(\U_n(x)-\widetilde \U_n(x)\right)>\varepsilon/2]
+ \limsup_{n\to\infty}\mathbb{P}[|c_n\left(\U_n(y)-\widetilde \U_n(y)\right)>\varepsilon/2].
\end{align*}
By  \eqref{eq:B2}, the right-hand side tends to zero as $x,y\to+\infty$. Thus, the family $(\intprocleft(x))_{x>1}$ satisfies the Cauchy criterion for convergence in probability. Denote its limit by $\intprocleft^0$.
	Theorem~3.2 in~\cite{Billingsley}, with approximation error controlled precisely
	by \eqref{eq:B2}, now gives \ref{ass:A3} with new fluctuation on the right hand side given by
	$\intprocleftup(z):=\intprocleft(z)-\intprocleft^0$ and which converges to $0$ in probability as $x\to\infty$.

\vspace{2mm}
\noindent	 
Conversely, assume \ref{ass:A3} and $\intprocleftup (x)\to0$ in probability, then we argued already in \eqref{eq:Psi_to_L} how \ref{ass:A2} follows. Moreover,  Portmanteau gives
$$
	\limsup_{n\to\infty}
	\mathbb P\!\left(
	c_n\left|\U_n(x)-\widetilde \U_n(x)\right|>\varepsilon
	\right)
	\le \mathbb{P}[|\intprocleftup (x)|\ge\varepsilon]\to 0\text{ as }x\to\infty.$$
\end{proof}
\section{Appendix}

\begin{proposition}[Vitali-Porter theorem]
Let $(f_n)_{n\in\mathbb{N}}$ be a locally bounded sequence of analytic functions defined on a domain $D\subseteq \mathbb{C}$ and such that the limit $\lim_{n\to\infty}f_n(z)$ exists for all $z\in E$, where $E\subset D$ is a subset having an accumulation point in $D$. Then $f_n$ converges locally uniformly on $D$ to an analytic function.
\end{proposition}

Here is a version of Vitali-Porter theorem for random analytic processes. Let $\mathcal{A}(D)$ denote the space of functions analytic on $D$ and endowed with the topology of locally uniform convergence.

\begin{theorem}\label{thm:Vitali-Porter-random}
Let $(Z_n)_{n\in\mathbb{N}}$ be a sequence of random elements taking values in $\mathcal{A}(D)$ and such that the corresponding sequence $(\mathbb{P}_n)_{n\in\mathbb{N}}$ of their distributions is tight. Let $E\subset D$ be a (non-random) subset of $D$ having an accumulation point in $D$. Assume that, for any $m\in\mathbb{N}$ and $t_1,\ldots,t_m\in E$, 
\begin{equation}\label{eq:vitali-porter_random_assumption1}
(Z_n(t_1),\ldots,Z_n(t_m))~\toweak~(Z_{\infty}(t_1),\ldots,Z_{\infty}(t_m)),
\end{equation}
for some collection of random variables $(Z_{\infty}(t))_{t\in E}$. Then,
\begin{itemize}
\item[(i)] the collection $(Z_{\infty}(t))_{t\in E}$ defines uniquely a distribution $\mathbb{P}_\infty$ on $\mathcal{A}(D)$;
\item[(ii)] $\mathbb{P}_n\toweak \mathbb{P}_\infty$ on the space of probability measures on $\mathcal{A}(D)$.
\end{itemize}
\end{theorem}
\begin{proof}
Recall that $\mathcal{A}(D)$, endowed with the topology of locally uniform convergence, is a Polish space. We first show that the finite-dimensional distributions on $E$ uniquely determine a probability measure on $\mathcal{A}(D)$. Choose a countable subset\footnote{If $E$ is itself countable we can choose $E_0=E$.} $E_0=\{t_1,t_2,\ldots\}\subset E$ having an accumulation point in $D$, and consider the map
$$
T:\mathcal{A}(D)\to\mathbb{C}^{\infty},\quad T(f)=(f(t_1),f(t_2),\ldots).
$$
The map $T$ is continuous with respect to the product topology on $\mathbb{C}^{\infty}$, and it is injective by the identity theorem for analytic functions. Since both $\mathcal{A}(D)$ and $\mathbb{C}^{\infty}$ are Polish spaces, the Lusin--Souslin theorem implies that $T(\mathcal{A}(D))$ is a Borel subset of $\mathbb{C}^{\infty}$ and that the inverse map $T^{-1}:T(\mathcal{A}(D))\to\mathcal{A}(D)$ is Borel measurable. Consequently,
$$
\sigma\bigl(\{f\mapsto f(t):t\in E_0\bigr)=\mathcal{B}(\mathcal{A}(D)),
$$
that is, the smallest $\sigma$-algebra making every evaluation map at points of $E_0$ measurable coincides with the Borel $\sigma$-algebra on $\mathcal{A}(D)$. Thus, any two probability measures on $\mathcal{A}(D)$ having the same finite-dimensional distributions on $E_0$ must coincide.

We now establish the existence of the limiting measure. By tightness and Pro\-kho\-rov's theorem, there exists a subsequence $(\mathbb{P}_{n_k})_{k\in\mathbb{N}}$ and a probability measure $\mathbb{P}_\infty$ on $\mathcal{A}(D)$ such that $\mathbb{P}_{n_k}\toweakk\mathbb{P}_\infty$. Moreover, the finite-dimensional projections of $\mathbb{P}_\infty$ on $E$ agree with the corresponding finite-dimensional distributions of $Z_{\infty}$ on $E$. By the uniqueness established above, it is the only probability measure on $\mathcal{A}(D)$ with this property. This proves (i).

To prove (ii), consider an arbitrary subsequence $(\mathbb{P}_{n_m})_{m\in\mathbb{N}}$. By tightness, it contains a further subsequence $(\mathbb{P}_{n_{m_k}})_{k\in\mathbb{N}}$ converging weakly to some probability measure $\mathbb{P}'_\infty$ on $\mathcal{A}(D)$. The continuous mapping theorem and~\eqref{eq:vitali-porter_random_assumption1} imply that $\mathbb{P}'_\infty$ has the same finite-dimensional distributions on $E$ as $\mathbb{P}_\infty$. By (i), we conclude that $\mathbb{P}'_\infty=\mathbb{P}_\infty$. Thus, we have proved that 
every subsequence of $(\mathbb{P}_n)$ contains a further subsequence converging weakly to $\mathbb{P}_\infty$, which means that the entire sequence converges. The proof is complete.
\end{proof}

\section*{Acknowledgments and AI-statement}
JJ has been financially supported by the DFG priority program SPP 2265 Random Geometric Systems. The authors are grateful to Zakhar Kabluchko for his valuable ideas and contributions during the early stages of this project, which played an important role in shaping the original approach to the problem. The authors obtained all main results before AI was available. Later on, the manuscript has been prompted to ChatGPT~6 for proof-checks and assistance in finalizing the exposition.

\bibliographystyle{plain}
\bibliography{symmetric}
\end{document}